\documentclass{article}

\usepackage{graphicx}
\usepackage[german,english]{babel}
\usepackage{amsmath}
\usepackage{amssymb}
\usepackage{amsfonts}
\usepackage{bbm}
\usepackage{url}
\usepackage{booktabs}
\usepackage{tabularx}
\usepackage{epsfig,subfigure}
\usepackage{longtable}
\usepackage{lscape}
\usepackage{psfrag}

\newtheorem{defn}{Definition}

\newtheorem{prop}[defn]{Proposition}
\newtheorem{thm}[defn]{Theorem}
\newtheorem{cor}[defn]{Corollary}
\newtheorem{conj}[defn]{Conjecture}
\newtheorem{ques}[defn]{Question}

\newtheorem{upthm}[defn]{Upper Bound Theorem}
\newtheorem{lowthm}[defn]{Lower Bound Theorem}
\newtheorem{mapthm}[defn]{Map Color Theorem}

\makeatletter

\graphicspath{{.}}

\title{Triangulated Manifolds with Few Vertices: Combinatorial Manifolds}

\author{\Large Frank H.~Lutz}
\date{}

\begin{document}

\selectlanguage{english}

\maketitle

\bigskip
\bigskip
\bigskip

\vfill

%%%%%%%%%%%%%%%%%%%%%%%%%%%%%%%%%%%%%%%%%%%%%%%%%%%%%%%%%%%%%%%%%%%

\noindent
There are essentially two ways to decompose 
a (compact, connected) $d$-mani\-fold (without boundary) into $d$-simplices:
as a (standard) simplicial complex or, more generally, as a simplicial cell
complex. In the latter case, identifications on the boundaries 
of the maximal simplices are sometimes allowed.

As an example, we need five tetrahedra to triangulate the $3$-sphere
as a simplicial complex, but only two tetrahedra to triangulate it
as a simplicial cell complex without identifications on the
boundaries of the tetrahedra. If we allow for identifications on the boundary,
then one tetrahedron suffices.

Triangulations of $3$-manifolds as general simplicial cell complexes
have been studied intensively in recent years (see for example 
\cite{MartelliPetronio2004}, \cite{Matveev1990}, \cite{Matveev1998}, \cite{Matveev2003}
and the references contained therein): 
The $3$-manifolds $S^3$, ${\mathbb R}{\bf P}^{\,3}$, and $L(3,1)$ can be triangulated
with one tetrahedron only. Otherwise the minimal number of tetrahedra
for a triangulation of a (closed) orientable irreducible $3$-manifold $M$
coincides with the Matveev complexity $c(M)$ of $M$; cf.\ \cite[p.~62]{Matveev2003}.

In the following, our focus will be on triangulations of
manifolds as proper simplicial complexes and on their combinatorial 
properties. A triangulated $d$-manifold $M$ is a \emph{simplicial $d$-manifold}
if it is triangulated as a simplicial complex. If, in addition, every of the
vertex-links of $M$ is PL homeomorphic to the boundary of a $d$-simplex,
then $M$ is a \emph{combinatorial (= triangulated PL) \linebreak 
$d$-manifold}.
%if it is a simplicial complex and if, in addition, every of its vertex-links is PL homeomorphic
%to the boundary of a $d$-simplex.

For combinatorial manifolds, various restrictions are known 
on the numbers of $i$-dimensional faces $f_i$,
in particular, on the number of vertices $n=f_0$.
In the first two sections of this 
% chapter,
paper, 
we give a survey on such restrictions.
Explicit examples of small (or otherwise interesting)  
triangulations of manifolds of dimension up to eight
will be discussed in the subsequent sections.

%\pagebreak

%\section{Lower Bounds and Minimal Triangulations}
\section{Minimal Triangulations}

Vertex-minimal triangulations are presently known only for a rather
limited number of manifolds due to the following:

\begin{itemize}
\item There are, apart from some special cases, 
      no standard constructions to produce triangulations of a given 
      manifold $M$ with \emph{few} vertices.\\[-2mm]
\item If such a small triangulation is obtained in some way with, say, $n$ vertices,
      then proving that $n$ is minimal is, in general, a difficult task.
      The only reasonable approach (apart from a complete enumeration)
      is to provide lower bounds on $n$ 
      in terms of numerical invariants of the given manifold~$M$. 
      Minimality follows if one of these bounds is sharp.
\end{itemize}
For $2$-dimensional manifolds $M$, Heawood \cite{Heawood1890} proved a lower bound in terms
of the Euler characteristic of $M$, and Jungerman and Ringel
settled the problem of constructing corresponding minimal triangulations.

Let a triangulation on $n$ vertices be \emph{$k$-neighborly} if 
every $k$-subset of the $n$ vertices forms a face of the
triangulation. In other words, the triangulation has the 
$(k-1)$-skeleton of the $(n-1)$-simplex on $n$ vertices.
We \underline{underline} the $k$-th entry of the \emph{face vector}
$f=(f_0,\dots,f_d)$ of a triangulation,
whenever the triangulation is $k$-neighborly.
\begin{thm} \label{thm:Heawood}
{\rm (Jungerman and Ringel \cite{JungermanRingel1980}, \cite{Ringel1955})}
If $M$ is $2$-dimensional and different from the orientable
surface of genus $2$, the Klein bottle, and the non-orient\-able surface of genus $3$,
then there is a triangulation of $M$ on\, $n$ vertices if and only if
\begin{equation}
\binom{n-3}{2}\geq 3\,(2-\chi (M)),
\end{equation}
with equality if and only if the triangulation is $2$-neighborly.
(For the three exceptional cases, $\binom{n-3}{2}$
has to be replaced by $\binom{n-4}{2}$, cf.~{\rm \cite{Huneke1978}}.)
\end{thm}
Otherwise the outcome is poor: 
Apart from minimal triangulations of the $d$-sphere as the boundary of the
$(d+1)$-simplex and one series of minimal triangulations of (twisted) sphere
products by K\"uhnel \cite{Kuehnel1986a-series} (see Theorem~\ref{thm:K-series}), 
there are merely eleven examples of manifolds
for which we have minimal triangulations
(with an explicit proof for minimality); see Table~\ref{tbl:minimal}.

Before we have a closer look at these examples in the following sections,
we first give an overview over all lower bounds that are currently known.  %in Table~\ref{tbl:minimal}, 
\begin{thm} \label{thm:bk-minimal} %\label{thm:minimal}
{\rm (Brehm and K\"uhnel~\cite{BrehmKuehnel1987})} Let $M$ be a combinatorial $d$-manifold with $n$ vertices.
\begin{itemize}
\item[{\rm (a)}] If $M$ is not a sphere, then
\begin{equation}
%n\geq 3\lceil\frac{d}{2}\rceil +3,
%n\geq 3\lceil\tfrac{d}{2}\rceil +3,
n\geq 3\left\lceil\frac{d}{2}\right\rceil +3,
\end{equation}
where equality can only occur for $d=2,4,8,16$. In these cases,
$M$ is a manifold `like a projective plane' (in the sense of\, {\rm\cite{EellsKuiper1962b}}).
\item[{\rm (b)}] If $M$ is $(i-1)$-connected but not $i$-connected and\, $1\leq i<d/2$, then
\begin{equation}
n\geq 2d+4-i.
\end{equation}
\end{itemize}
\end{thm}
There is a combinatorially unique triangulation ${\mathbb R}{\bf P}^{\,2}_6$ of the real projective plane 
with $6$ vertices and a combinatorially unique triangulation ${\mathbb C}{\bf P}^{\,2}_9$  
of the complex projective plane with $9$ vertices \cite{KuehnelBanchoff1983}, \cite{KuehnelLassmann1983-unique}.
Moreover, there are at least six combinatorially different $15$-vertex triangulations 
of an $8$-dimensional mani\-fold ${\sim}{\mathbb H}{\bf P}^{\,2}$ 
like the quaternionic projective plane; for details see
Sections~\ref{sec:2-manifolds}, \ref{sec:4-manifolds}, and \ref{sec:8-manifolds}.
\begin{ques} {\rm (Brehm and K\"uhnel~\cite{BrehmKuehnel1987})} \index{projective space!Cayley}
Can the $16$-dimensional Cayley projective plane
be triangulated with $27$ vertices?
\end{ques}

If $d$ is even and $M$ is not $d/2$-connected, then Part~(a) 
of Theorem~\ref{thm:bk-minimal} yields
$n\geq 3d/2+3$ for manifolds `like a projective plane' 
and $n\geq 3d/2+4$ for all other manifolds.

\begin{cor} {\rm (Brehm and K\"uhnel~\cite{BrehmKuehnel1987})} \label{cor:BrehmKuehnel}
Let $M$ be a combinatorial $d$-manifold with $n$ vertices.
\begin{itemize}
\item[{\rm (a)}] If $M$ has the same homology as\, $S^{d-i}\!\times\!S^i$, then $n\geq 2d+4-i$.\\[-2.5mm]
\item[{\rm (b)}] If $M$ is not simply connected, then $n\geq 6$ for $d=2$ 
                 and $n\geq 2d+3$ for $d\geq 3$.
\end{itemize}
\end{cor}

Theorem~\ref{thm:bk-minimal}(b) proves vertex-minimality for
a series of triangulations of (twisted) sphere products
constructed by K\"uhnel.

\begin{thm} {\rm (K\"uhnel \cite{Kuehnel1986a-series})}\label{thm:K-series}
There is a series $M^d$ of vertex-minimal combinatorial triangulations 
with transitive dihedral auto\-morphism group on\, $2d+3$\, vertices
of topological type
\begin{itemize}
\item $M^d\cong S^{d-1}\!\times\!S^1$, if\, $d$\, is even, and\\[-2.5mm]
\item $M^d\cong S^{d-1}\hbox{$\times\hspace{-1.62ex}\_\hspace{-.4ex}\_\hspace{.7ex}$}S^1$, if\, $d$\, is odd.
\end{itemize}
\end{thm}
This series was later generalized by K\"uhnel and Lassmann \cite{KuehnelLassmann1996-bundle}
to a series of combinatorial $d$-manifolds $M^d_k(n)$, $1\leq k\leq d-1$, 
with vertex-transitive dihedral symmetry group on $n\geq 2^{d-k}(k+3)-1$
vertices. In general, $M^d_k(n)$ is a (orientable or non-orientable) $k$-dimensional sphere
bundle over the $(d-k)$-dimensional torus with $M^d_{d-1}(2d+3)=M^d$
and $M^d_{d-1}(2d+4)\cong S^{d-1}\!\times\!S^1$ as special cases.
\begin{conj}\label{conj:products_even_dim}
The transitive K\"uhnel-Lassmann triangulation $M^d_{d-1}(2d+4)$
of\,\, $S^{d-1}\!\times\!S^1$\, on\, $2d+4$\, vertices is
vertex-minimal for odd $d$. 
For even $d$, $S^{d-1}\hbox{$\times\hspace{-1.62ex}\_\hspace{-.4ex}\_\hspace{.7ex}$}S^1$\,
can be triangulated minimally  on\, $2d+4$\, vertices. 
\end{conj}
The Klein bottle needs at least $8$ vertices for a triangulation \cite{Franklin1934},
so the conjecture holds for $d=2$. For $d=3$ it follows from Theorem~\ref{thm:minimalthree}
of Walkup below. Moreover, $S^{3-1}\hbox{$\times\hspace{-1.62ex}\_\hspace{-.4ex}\_\hspace{.7ex}$}S^1$
has triangulations with $12$ vertices, but minimality is not settled in
this case; see Section~\ref{sec:4-manifolds}.

\begin{prop}\label{prop:series_2d+3}
There are triangulations of\, $S^{d-1}\hbox{$\times\hspace{-1.62ex}\_\hspace{-.4ex}\_\hspace{.7ex}$}S^1$\,
with\, $3d+3$\, vertices.
\end{prop}

\noindent
\textbf{Proof.} 
%\begin{proof}
Let $I$ be an interval, triangulated with $4$ vertices.
If we take the boundary of the $d$-simplex with $d+1$ vertices as a 
triangulation of the sphere $S^{d-1}$, then the product triangulation
(see 
%Chapter~\ref{ch:geometric}) 
\cite{Lutz2003bpre})
of\, $S^{d-1}\times I$ has\, $4(d+1)$ vertices. 
The boundary of $S^{d-1}\!\times I$ consists of two disjoint copies 
of the boundary of the $d$-simplex. By identifying the two boundary spheres 
(with taking care of the orientation), we obtain the wanted triangulation of 
$S^{d-1}\hbox{$\times\hspace{-1.62ex}\_\hspace{-.4ex}\_\hspace{.7ex}$}S^1$
with $3d+3$ vertices. 
%\end{proof}
\hfill $\Box$
\bigskip

The Brehm-K\"uhnel bound of Theorem~\ref{thm:bk-minimal}(a) is sharp for $d=3,4,8$: 
$S^2\hbox{$\times\hspace{-1.62ex}\_\hspace{-.4ex}\_\hspace{.7ex}$}S^1$
can be triangulated with $9$ \cite{Walkup1970b}, ${\mathbb C}{\bf P}^{\,2}$ with $9$ \cite{KuehnelBanchoff1983}, 
and ${\sim}{\mathbb H}{\bf P}^{\,2}$ with $15$ vertices \cite{BrehmKuehnel1992}. 
In Sections~\ref{sec:5-manifolds} and \ref{sec:6-manifolds} we give triangulations
of $S^3\!\times\!S^2$ with $12$ and of $S^3\!\times\!S^3$ with $13$
vertices, and therefore show that the Brehm-K\"uhnel bound
is also sharp in dimensions $5$ and $6$.
\begin{ques}
Is there a $15$-vertex triangulation of $S^4\!\times\!S^3$?
\end{ques}
Such a triangulation would settle the tightness of the Brehm-K\"uhnel
bound for $d=7$; it would be a minimal triangulation also
by Corollary~\ref{cor:BrehmKuehnel}(a).
%and it would be a further example of a tight triangulation 
%(see Chapter~\ref{ch:tight}).

Starting from a product triangulation with $6\cdot 4=24$ vertices and
by applying bistellar flips to it (cf.\ \cite{BjoernerLutz2000}
for a discussion of bistellar flips), \linebreak
we found a triangulation of $S^4\!\times\!S^2$ with $16$ vertices.
\begin{ques}
Can $S^4\!\times\!S^2$ be triangulated with $14$ vertices?
\end{ques}

For homology spheres, Bagchi and Datta improved Brehm and K\"uhnel's bounds
of Theorem~\ref{thm:bk-minimal} in dimensions $3\leq d\leq 6$.

\begin{thm} Let $M$ be a combinatorial $d$-manifold with $n$ vertices.\label{thm:homology_sphere}
\begin{itemize}
\item[{\rm (a)}] {\rm (Brehm and K\"uhnel~\cite{BrehmKuehnel1987})} 
                 If $M$ is a ${\mathbb Z}$-homology sphere, then $n\geq 2d+3$ for\, $d\geq 6$.\\[-1mm]
\item[{\rm (b)}] {\rm (Bagchi and Datta~\cite{BagchiDatta2003pre})} 
                 If $M$ is a ${\mathbb Z}_2$-homology sphere, then $n\geq d+9$ for\, $3\leq d\leq 6$.
                 If $n=d+9$, then $M$ does not admit a (non-trivial) bistellar flip.
\end{itemize}
\end{thm}

The result of Bagchi and Datta implies that at least $12$ vertices
are needed to triangulate the lens space $L(3,1)$. A triangulation
with this number of vertices was first found by Brehm~\cite{Brehm-pers};
see Section~\ref{sec:3-manifolds}.

\pagebreak

The minimal number of vertices is known for four further $3$-manifolds.

\begin{thm}\label{thm:minimalthree}\mbox{}
{\rm (Walkup~\cite{Walkup1970})} With the exception of $S^3$,
$S^2\hbox{$\times\hspace{-1.62ex}\_\hspace{-.4ex}\_\hspace{.7ex}$}S^1$, 
and \hbox{$S^2\!\times\!S^1$}, which have minimal
triangulations with $5$, $9$ and $10$ vertices, respectively, every
other triangulated $3$-manifold has at least $11$ vertices.
% The real projective space 
${\mathbb R}{\bf P}^{\,3}$ can be triangulated with
$11$ vertices.
\end{thm}

For real and complex projective spaces there are lower bounds due to
Arnoux and Marin.

\begin{thm} {\rm (Arnoux and Marin~\cite{ArnouxMarin1991})} \label{thm:ArnouxMarin}
Combinatorial triangulations of\, ${\mathbb R}{\bf P}^d$\, and\, ${\mathbb C}{\bf P}^r$
have at least\, $n\geq (d+1)(d+2)/2$\, and\, $n\geq (r+1)^2$\, 
%have at least\, $n\geq {d+2\choose 2}$\, and\, $n\geq (r+1)^2$\, 
vertices respectively. Equality is possible only for $d=2$ and $r=2$.
\end{thm}

%It therefore also follows from Arnoux and Marin's bound that
%Walkup's triangulation of ${\mathbb R}{\bf P}^{\,3}$ with $11$ vertices
%is vertex-minimal.

A series of combinatorial triangulations of ${\mathbb R}{\bf P}^{\,d}$ 
with $2^{\,d+1}-1$ vertices was constructed by K\"uhnel~\cite{Kuehnel1986b-kummer}.
We applied the bistellar flip program BISTELLAR~\cite{Lutz_BISTELLAR} to the $4$-dimensional example and obtained
a $16$-vertex triangulation of ${\mathbb R}{\bf P}^{\,4}$ (see 
%Chapter~\ref{ch:pseudo}),
\cite{Lutz2004dpre}),
which is minimal by Theorem~\ref{thm:ArnouxMarin}.
Furthermore, we found a triangulation of ${\mathbb R}{\bf P}^{\,5}$
with $24$ vertices.
\begin{ques}
Is there a $22$-vertex triangulation of\, ${\mathbb R}{\bf P}^{\,5}$?
\end{ques}
For the complex projective spaces ${\mathbb C}{\bf P}^r$ 
explicit triangulations are not known for $r>2$.
\begin{ques}
Is there a $17$-vertex triangulation of\, ${\mathbb C}{\bf P}^{\,3}$?
\end{ques}

A $4$-dimensional analogue of Heawood's bound for surfaces from Theorem~\ref{thm:Heawood}
was proved by K\"uhnel.

\begin{thm} {\rm (K\"uhnel \cite[4.1]{Kuehnel1990-few})} \label{thm:kuhnel-4dim}
If $M$ is a combinatorial $4$-manifold with $n$ vertices, then
\begin{equation}\label{eq:bound4}
\binom{n-4}{3}\geq 10\,(\chi (M)-2).
\end{equation}
Equality holds if and only if $M$ is $3$-neighborly.
\end{thm}

In generalization of Theorem~\ref{thm:Heawood} and Theorem~\ref{thm:kuhnel-4dim}, 
K\"uhnel and, independently, Kalai conjectured:
\begin{conj} {\rm \cite[p.~61]{Kuehnel1990-few}}\label{conj:kuehnel_kalai}
Let $M$ be a combinatorial $2k$-manifold with $n$ vertices, then
\begin{equation}\label{eq:bound}
\displaystyle \binom{n-k-2}{k+1}\geq\,(-1)^k \binom{2k+1}{k+1} \big(
\chi (M) - 2 \big),
\end{equation}
with equality if and only if the triangulation is $(k+1)$-neighborly.
\end{conj}
The conjecture holds in several cases (cf.~\cite[4.8]{Kuehnel1995-book}), 
in particular, for $n\leq 3k+3$ and $n\geq k^2+4k+3$.
The latter inequality was improved to $n\geq 4k+3$ by 
Novik~\cite{Novik1998}. Moreover, Novik showed that for both ranges
the conjecture holds for arbitrary simplicial triangulations, not just combinatorial ones(!). 
As the bounds~\ref{eq:bound4} and \ref{eq:bound} depend on the Euler characteristic,
we get more detailed information on the minimal numbers of vertices $n$
for triangulations of $2k$-man\-i\-folds than by the general bounds
of Theorem~\ref{thm:bk-minimal}.
Since the Euler characteristic is zero for
odd-di\-men\-sion\-al manifolds, it was proposed by Kalai
\cite{BjoernerKalai1989}, \cite{Kalai1987}
to replace the factor $(-1)^k (\chi (M) - 2)$ by the sum
of the Betti numbers $\sum_{i=1}^{2k-i}\beta_i(M)$ (or possibly by a weighted sum). 
By Poincar\'e duality,
$(-1)^k(\chi -2)=\beta_k-(\beta_{k+1}+\beta_{k-1})+(\beta_{k+2}+\beta_{k-2})-\cdots =\beta_k+2\sum_{i=1}^{k-1}(-1)^i\beta_{k-i}$\,
and\, $\sum_{i=1}^{2k-1}\beta_i=\beta_k+2\sum_{i=1}^{k-1}\beta_{k-i}$,
when the Betti numbers $\beta_i$ are computed with respect to the 
field ${\mathbb F}_2$ ($M$ is orientable over ${\mathbb F}_2$!).

\medskip

\begin{thm} {\rm (Novik \cite{Novik1998})}
Let $M$ be a simplicial $d$-manifold with $n$ vertices.
\begin{itemize}
\item[{\rm (1)}]
If $d=2k$ and\, $n\leq 3k+3$\, or\, $n\geq 4k+3$, then
\begin{equation}
\displaystyle \binom{n-k-2}{k+1}\geq\,\binom{2k+1}{k+1} \big( \beta_k+2\sum_{i=0}^{k-2}\tilde{\beta}_{i}\big).
\end{equation}
\item[{\rm (2)}]
If $d=2k$ and\, $n\leq 3k+3$\, or\, $n\geq 7k+3$, then
\begin{equation}
\displaystyle \binom{n-k-2}{k+1}\geq\,\binom{2k+1}{k+1} \big( \beta_k+2\sum_{i=1}^{k-1}\beta_{i}\big).
\end{equation}
\item[{\rm (3)}]
If $d=2k-1$ and\, $n\leq 3k+2$\, or\, $n\geq 4k+1$, then
\begin{equation}
\displaystyle \frac{2n}{n+k+2}\cdot\binom{n-k-2}{k}\geq\,\binom{2k-1}{k} \big( 2\sum_{i=1}^{k-1}\beta_{i}\big).
\end{equation}
\end{itemize}
\end{thm}

K\"uhnel \cite{Kuehnel-pers} has worked out
another elegant (and previously unpublished) conjecture 
that generalizes Theorem~\ref{thm:bk-minimal}(b) (in the case of (twisted) sphere products) and
Conjecture~\ref{conj:kuehnel_kalai} (in the case of $(k-1)$-connected $2k$-manifolds).
\begin{conj}  {\rm (K\"uhnel)} \label{conj:kuehnel_pascal}
Let $M$ be a combinatorial $d$-manifold with $n$ vertices
and reduced Betti numbers $\tilde{\beta}_j(M)$, $j=0,\dots,\lfloor\frac{d}{2}\rfloor$. Then
\begin{equation}\label{eq:bound_k1}
\displaystyle \binom{n-d+j-2}{j+1}\geq\, \binom{d+2}{j+1}\,\tilde{\beta}_{j}(M)\quad\mbox{for}\quad\textstyle j=0,\dots,\lfloor\frac{d-1}{2}\rfloor.
\end{equation}\\[-1mm]
If $d$ is even, then additionally (for\, $j=\frac{d}{2}$):
\begin{equation}\label{eq:bound_k2}
\displaystyle \binom{n-\frac{d}{2}-2}{\frac{d}{2}+1}\geq\,
\binom{d+2}{\frac{d}{2}+1}\,\frac{\tilde{\beta}_{\frac{d}{2}}(M)}{2}.
\end{equation}
If equality holds for one of the bounds of~\ref{eq:bound_k1}
or for~\ref{eq:bound_k2}, say, for $j=s$, then $\tilde{\beta}_j(M)=0$
for $j\neq s$.
\end{conj}
The bounds of K\"uhnel can be arranged in a Pascal-like triangle
as displayed in Table~\ref{tbl:conj_kuehnel}: There is a row 
for every dimension $d=0,1,\dots$ with the respective bounds 
for $j=0,\dots,\lceil\frac{d-1}{2}\rceil$
as entries.

\begin{table}
\small\centering
\defaultaddspace=0.15em
\caption{The Pascal-like triangle of K\"uhnel's lower
  bounds of Conjecture~\ref{conj:kuehnel_pascal}.}\label{tbl:conj_kuehnel}
\begin{tabular*}{\linewidth}{@{\extracolsep{\fill}}l@{\hspace{1mm}}l@{\hspace{1mm}}l@{\hspace{1mm}}l@{\hspace{1mm}}l@{\hspace{1mm}}l@{\hspace{1mm}}l@{\hspace{1mm}}l}
 \addlinespace
 \addlinespace
 \addlinespace
 \addlinespace
\toprule
 \addlinespace
 \addlinespace
 \addlinespace
 \addlinespace
 \addlinespace
                        &                                  & &                        & & $\binom{n-2}{1}\geq 0$ \\
 \addlinespace
                        &                                  & &                & $\binom{n-3}{1}\geq 0$ \\
 \addlinespace
                        &                                  & & $\binom{n-4}{1}\geq 0$ & & $\binom{n-3}{2}\geq\binom{4}{2}\frac{\beta_1}{2}$ \\
 \addlinespace
                        &         & $\binom{n-5}{1}\geq 0$ & & $\binom{n-4}{2}\geq\binom{5}{2}\beta_1$ \\
 \addlinespace
                        & $\binom{n-6}{1}\geq 0$ & & $\binom{n-5}{2}\geq\binom{6}{2}\beta_1$ & & $\binom{n-4}{3}\geq\binom{6}{3}\frac{\beta_2}{2}$ \\
 \addlinespace
             $\binom{n-7}{1}\geq 0$ & & $\binom{n-6}{2}\geq\binom{7}{2}\beta_1$ & & $\binom{n-5}{3}\geq\binom{7}{3}\beta_2$ \\
 \addlinespace
  & \hspace{4mm}\dots & & \hspace{4mm}\dots & & \hspace{4mm}\dots \\

% \addlinespace
%6 & $\binom{n-8}{1}\geq 0$ & & $\binom{n-7}{2}\geq\binom{8}{2}\beta_1$ & & $\binom{n-6}{3}\geq\binom{8}{3}\beta_2$ & & $\binom{n-5}{4}\geq\binom{8}{4}\frac{\beta_3}{2}$ \\
 \addlinespace
 \addlinespace
 \addlinespace
\bottomrule
\end{tabular*}
\end{table}

The diagonal in the K\"uhnel triangle corresponding to $j=0$ gives the
bound $\binom{n-d-2}{1}\geq 0$,
which is equivalent to $n\geq d+2$. In other words, every
combinatorial $d$-manifold trivially has at least as many vertices as
the boundary of a $(d+1)$-simplex.

For $d$-manifolds $M$ with $\beta_1(M)=1$, 
the second diagonal $j=1$ in  K\"uhnel's triangle 
yields (for $d\geq 3$) the bound $\binom{n-d-1}{2}\geq\, \binom{d+2}{2}$, or, equivalently,
\mbox{$n\geq 2d+3$}. This is precisely the bound for non-simply connected 
$d$-manifolds from Corollary~\ref{cor:BrehmKuehnel}(b); K\"uhnel's
series from Theorem~\ref{thm:K-series} of vertex-minimal
triangulations of (twisted) $S^{d-1}$-bundles over $S^1$ 
with $n= 2d+3$ vertices shows that this bound is best possible.

In fact, the bounds of Conjecture~\ref{conj:kuehnel_pascal} hold
for all sphere products \mbox{$S^{d-i}\!\times\! S^i$}, for which the relevant
bound $\binom{n-d+i-2}{i+1}\geq \binom{d+2}{i+1}$
gives $n\geq 2d+4-i$ in accordance with Corollary~\ref{cor:BrehmKuehnel}(a).
For $(k-1)$-connected $2k$-manifolds the bounds of Conjecture~\ref{conj:kuehnel_pascal}
corresponding to $j=k$ coincide with the generalized Heawood bounds 
of Conjecture~\ref{conj:kuehnel_kalai}.

\emph{Let us point out that K\"uhnel's Conjecture~\ref{conj:kuehnel_pascal} is of particular
interest for odd-dimensional manifolds $M$,
where we, up to now, only have the bounds of Theorem~\ref{thm:bk-minimal}
(together with the special bounds of Theorem~\ref{thm:homology_sphere}(b) 
and of Theorem~\ref{thm:ArnouxMarin}).}

A first case with $\beta_1(M)>1$, where K\"uhnel's Conjecture~\ref{conj:kuehnel_pascal}
would prove vertex-mini\-ma\-li\-ty, is for the manifolds
$(S^2\!\times\!S^1)^{\#3}$ and\, $(S^2\hbox{$\times\hspace{-1.62ex}\_\hspace{-.4ex}\_\hspace{.7ex}$}S^1)^{\#3}$
for which we know triangulations with $13$ vertices.
For $4$-manifolds $M$ with trivial second homo\-logy group $H_2(M)=0$,
K\"uhnel's bound $\binom{n-5}{2}\geq 15\cdot\beta_1(M)$ follows from
results of Walkup (cf.\ \cite[p.~96]{Kuehnel1995-book}, \cite[Thm.~5]{Walkup1970b}).

As mentioned above, there are, besides the minimal triangulations by Jungerman and Ringel
of surfaces, the $d$-sphere, which can be
triangulated as the boundary of the $(d+1)$-simplex,
and the series of minimal triangulations of (twisted) sphere products 
by K\"uhnel \cite{Kuehnel1986a-series}, eleven exceptional examples of
manifolds for which we explicitly have minimal triangulations; see Table~\ref{tbl:minimal}.

\pagebreak

\noindent
Minimality of these examples follows
\begin{itemize}
\item for\, $S^2\!\times\!S^1$\, and ${\mathbb R}{\bf P}^{\,3}=L(2,1)$\, from Theorem~\ref{thm:minimalthree},\\[-2.5mm]
\item for\, $L(3,1)$\, from Theorem~\ref{thm:homology_sphere}(b),\\[-2.5mm]
\item for\, ${\mathbb C}{\bf P}^{\,2}$\, from Theorem~\ref{thm:bk-minimal}(a) as well as from Theorem~\ref{thm:kuhnel-4dim},\\[-2.5mm]
\item for\, $S^2\!\times\!S^2$\, from \cite{KuehnelLassmann1983-unique},\\[-2.5mm]
\item for\, $(S^2\!\times\!S^2)\# (S^2\!\times\!S^2)$\, and the\, K3 surface\, from Theorem~\ref{thm:kuhnel-4dim},\\[-2.5mm]
\item for\, ${\mathbb R}{\bf P}^{\,4}$\, from Theorem~\ref{thm:ArnouxMarin},\\[-2.5mm]
\item for\, $S^3\!\times\!S^2$\, from Theorem~\ref{thm:bk-minimal}(a) as well as from (b), and\\[-2.5mm]
\item for\, $S^3\!\times\!S^3$\, and\, ${\sim}{\mathbb H}{\bf P}^{\,2}$\, from Theorem~\ref{thm:bk-minimal}(a).
\end{itemize}

\begin{table}
\small\centering
\defaultaddspace=0.1em
\caption{Eleven exceptional cases of vertex-minimal triangulations.}\label{tbl:minimal}
%\caption{\protect\parbox[t]{8cm}{Eleven exceptional examples of vertex-minimal triangulations}
%                                 in dimensions\, $3\leq d\leq 8$.}\label{tbl:minimal}
\begin{tabular*}{\linewidth}{@{\extracolsep{\fill}}clrl@{}}
 \addlinespace
 \addlinespace
 \addlinespace
 \addlinespace
 \addlinespace
 \addlinespace
\toprule
 \addlinespace
 \addlinespace
 \addlinespace
Dimension &  Manifold          &        $n$       &  Reference \\ \midrule
 \addlinespace
 \addlinespace
 \addlinespace
    3     %& $S^2\hbox{$\times\hspace{-1.62ex}\_\hspace{-.4ex}\_\hspace{.7ex}$}S^1$ & $9$ & \cite[$N^9_{51}$]{AltshulerSteinberg1974}, \cite{Walkup1970} \\
% \addlinespace
          & $S^2\!\times\!S^1$  & $10$ & \cite{Walkup1970} \\
 \addlinespace
          & ${\mathbb R}{\bf P}^{\,3}=L(2,1)$ & $11$ & \cite[$S_{2\cdot 4}$]{BrehmSwiatkowski1993}, \cite{Walkup1970} \\
 \addlinespace
          & $L(3,1)$            & $12$ & \cite{BagchiDatta2003pre}, \cite{Brehm-pers} \\
 \addlinespace
    4     & ${\mathbb C}{\bf P}^{\,2}$ & $9$ & \cite{KuehnelBanchoff1983}, \cite{KuehnelLassmann1983-unique} \\
 \addlinespace
%          & $S^3\!\times\!S^1$  & $11$ & \cite{Kuehnel1986a-series} \\
% \addlinespace
          & $S^2\!\times\!S^2$  & $11$ & \\%Ch.~\ref{ch:centrally} \\
 \addlinespace
          & $(S^2\!\times\!S^2)\# (S^2\!\times\!S^2)$ & $12$ & \\%Ch.~\ref{ch:manifold}\\
 \addlinespace
          & ${\mathbb R}{\bf P}^{\,4}$ & $16$ & \\%Ch.~\ref{ch:pseudo} \\
 \addlinespace
          & K3 surface          & $16$ & \cite{CasellaKuehnel2001} \\
 \addlinespace
    5     & $S^3\!\times\!S^2$  & $12$ & \\%Ch.~\ref{ch:centrally} \\
 \addlinespace
%          & $S^4\hbox{$\times\hspace{-1.62ex}\_\hspace{-.4ex}\_\hspace{.7ex}$}S^1$ & $13$ & \cite{Kuehnel1986a-series} \\
% \addlinespace
    6     & $S^3\!\times\!S^3$  & $13$ & \\%Ch.~\ref{ch:centrally} \\
 \addlinespace
%          & $S^5\!\times\!S^1$  & $15$ & \cite{Kuehnel1986a-series} \\
% \addlinespace
%    7     & $S^6\hbox{$\times\hspace{-1.62ex}\_\hspace{-.4ex}\_\hspace{.7ex}$}S^1$ & $17$ & \cite{Kuehnel1986a-series} \\
% \addlinespace
    8     & ${\sim}{\mathbb H}{\bf P}^{\,2}$ & $15$ & \cite{BrehmKuehnel1992} \\
 \addlinespace
%          & $S^7\!\times\!S^1$  & $19$ & \cite{Kuehnel1986a-series} \\
% \addlinespace
 \addlinespace
 \addlinespace
 \addlinespace
\bottomrule
\end{tabular*}
\end{table}

For various other small triangulations minimality cannot be proven 
due to a lack of sharper lower bounds.
\begin{conj}\label{conj:mintri}
The minimal number of vertices is
\begin{itemize}
\item[{\rm (a)}] $12$\, for\, $(S^2\!\times\!S^1)^{\#2}$ and\, $(S^2\hbox{$\times\hspace{-1.62ex}\_\hspace{-.4ex}\_\hspace{.7ex}$}S^1)^{\#2}$,\\[-2mm]
\item[{\rm (b)}] $13$\, for\, $(S^2\!\times\!S^1)^{\#3}$ and\, $(S^2\hbox{$\times\hspace{-1.62ex}\_\hspace{-.4ex}\_\hspace{.7ex}$}S^1)^{\#3}$,\\[-2mm]
\item[{\rm (c)}] $14$\, for the lens spaces\, $L(4,1)$, $L(5,2)$, and\, ${\mathbb R}{\bf P}^{\,2}\!\times\!S^1$,\\[-2mm]
\item[{\rm (d)}] $15$\, for\, $L(5,1)$, the prism spaces\, $P(2)=S^3/Q$, $P(3)$, $P(4)$,\,
                  and the connected sum\, ${\mathbb R}{\bf P}^{\,3}\#\,{\mathbb R}{\bf P}^{\,3}$, \\[-2mm]
\item[{\rm (e)}] $12$\, for\,\, ${\mathbb C}{\bf P}^{\,2}\#\,{\mathbb C}{\bf P}^{\,2}$\, and\,\,\, ${\mathbb C}{\bf P}^{\,2}\# -{\mathbb C}{\bf P}^{\,2}$,\\[-2mm]
\item[{\rm (f)}] $15$\, for\, $(S^3\hbox{$\times\hspace{-1.62ex}\_\hspace{-.4ex}\_\hspace{.7ex}$}S^1)\#\,({\mathbb C}{\bf P}^{\,2})^{\# 5}$, and\\[-2mm]
\item[{\rm (g)}] $13$\, for\,\, $SU(3)/SO(3)$.
\end{itemize}
\end{conj}
Triangulations for these examples with the given numbers of vertices are discussed in the
following sections.

\begin{thm} {\rm (K\"uhnel and Lassmann \cite{KuehnelLassmann1988-dtori})}
The $d$-dimensional torus $T^d$ can be triangulated combinatorially
with $2^{d+1}-1$ vertices (for $d\geq 2$).
\end{thm}
In fact, the corresponding series is a subseries of the K\"uhnel and Lassmann
series $M^d_k(n)$ with\, $T^d\cong M^d_1(2^{d+1}-1)$.

We tried to improve on the respective triangulations of $T^3$ and $T^4$
with our bistellar flip program, but with no success.
\begin{conj}\label{conj:Td}
The K\"uhnel-Lassmann series of combinatorial $d$-tori $T^d$ with\, $2^{d+1}-1$ vertices
is vertex-minimal.
\end{conj}
Besides the K\"uhnel-Lassmann triangulation of the $4$-torus,
a non-equivalent triangulation of $T^4$ with $31$ vertices was constructed 
by Grigis \cite{Grigis1998}.
Dartois and Grigis gave $25$ further triangulations of the eight-dimensional 
torus $T^8$ with $511$ vertices \cite{DartoisGrigis2000},
but so far, no additional triangulation of $T^3$ has been found with $15$ vertices.
\begin{conj}
The K\"uhnel-Lassmann triangulation {\rm \cite{KuehnelLassmann1984-3torus}} of $T^3$ with $15$ vertices 
is vertex-minimal and unique with $f=(15,\underline{105},180,90)$.
\end{conj}

\mathversion{bold}
\section{$f$-Vectors}
\mathversion{normal}

Of course, it would be desirable to know a complete characterization of the
$f$-vectors of a triangulable manifold, not just the minimal number of vertices.

For various classes of simplicial complexes this rather general question 
has, in fact, been answered (cf.\ the surveys by Bj\"orner \cite{Bjoerner1996} and Billera and
Bj\"orner~\cite{BilleraBjoerner1997}). For the class of all simplicial complexes
the result is known as the Kruskal-Katona theorem (\cite{Katona1968}, \cite{Kruskal1963};
see also Sch\"utzenberger \cite{Schuetzenberger1959}, who gave the first proof).
A complete characterization for simplicial polytopes 
is given by the $g$-theorem (conjectured by McMullen~\cite{McMullen1971})
of Billera and Lee \cite{BilleraLee1981} and Stanley \cite{Stanley1975}.
It is, however, an open problem whether the $g$-theorem can be extended to general
simplicial spheres (`$g$-conjecture for simplicial spheres' by McMullen \cite{McMullen1971}).
The $g$-conjecture holds for $3$- and $4$-dimensional spheres;
cf.\ McMullen \cite{McMullen1971}: There, the essential condition is that
$f_1\geq (d+1)n-\binom{d+2}{2}$, which is the case $k=1$ in the lower bound
theorem of Kalai and Gromov.
\begin{lowthm} {\rm (Kalai and Gromov, cf.\ \cite{Kalai1987} including
    the note added in proof and \cite[Ch.\ 2.4.10]{Gromov1986})}\index{lower bound theorem}
If $M$ is a\, $d$-di\-men\-sion\-al \mbox{(pseudo-)} manifold, then for every
triangulation with $n$ vertices,
\begin{equation}
f_k\,\geq\,\left\{
\begin{array}{l@{\hspace{3.5mm}}l@{\hspace{3.5mm}}l}
\,\binom{d+1}{k}n-\binom{d+2}{k+1}k & for & 1\leq k\leq d-1,\\[2.5mm]
\,d\cdot n-(d-1)(d+2)               & for & k=d.
\end{array}
\right.
\end{equation}
\end{lowthm}
This result for polytopal spheres was obtained before by Barnette~\cite{Barnette1973a},
the case\, $k=d$\, for pseudomanifolds by Klee~\cite{Klee1975}.
The lower bound theorem for $3$- and $4$-dimensional manifolds
was first proved by Walkup \cite{Walkup1970}.

\begin{upthm} {\rm (Novik \cite{Novik1998})}\index{upper bound theorem}
Let $M$ be a simplicial (homo\-logy) $d$-manifold.
If $d$ is odd and also if $d=2k$ and 
$\beta_k\leq 2\beta_{k-1}+2\sum_{i=1}^{k-3}\beta_k$,
then $M$ cannot have more $i$-dimensional faces as the corresponding cyclic
polytope with the same number of vertices, 
that is,\, $f_k(M)\leq f_k(C_{d+1}(n))$\, for\,\, $1\leq k\leq d$.

\end{upthm}
The upper bound theorem for polytopal spheres was first proved 
by McMullen \cite{McMullen1970} and later generalized to simplicial
spheres by Stanley~\cite{Stanley1975}.

%%%%%%%%%%%%%%%%%%%%%%%%%%%%%%%%%%%%%%%%%%%%%%%%%%%%%%%
%%%%%%%%%%%%%%%%%%%%%%%%%%%%%%%%%%%%%%%%%%%%%%%%%%%%%%%

\section{2-Manifolds}
\label{sec:2-manifolds}

The $f$-vectors of triangulated surfaces $M$ have a particularly
simple description.
We have Euler's equation, $n-f_1+f_2=\chi (M)$, and, by 
double counting of incidences between edges and triangles, 
it follows that $2f_1=3f_2$. Thus, the number of vertices 
$n$ determines $f_1$ and $f_2$, that is, a triangulated $2$-mani\-fold~$M$ 
of Euler characteristic $\chi (M)$
on $n$ vertices has $f$-vector
\begin{equation}
f=(n,3n-3\chi(M),2n-2\chi(M)).
\end{equation}
The \emph{orientable surface} $M(g,+)$ of genus $g$ has
homology 
\begin{equation}
H_*(M(g,+))=({\mathbb Z},{\mathbb Z}^{2g},{\mathbb Z})
\end{equation}
and Euler characteristic $\chi(M(g,+))=2-2g$,
whereas the \emph{non-orientable surface} $M(g,-)$ of genus $g$ has
homology 
\begin{equation}
H_*(M(g,-))=({\mathbb Z},{\mathbb Z}^{g-1}\oplus{\mathbb Z}_2,0)
\end{equation}
and Euler characteristic $\chi(M(g,-))=2-g$. 
The smallest possible $n$ for a triangulation of a $2$-manifold $M$
(with the exception of  the orientable
surface of genus~$2$, the Klein bottle, and the non-orient\-able
surface of genus $3$, where one more vertex is needed)
is determined by Heawood's bound 
\begin{equation}
n\geq\frac{1}{2}(7+\sqrt{49-24\chi (M)})
\end{equation}
of Theorem~\ref{thm:Heawood}. Corresponding minimal triangulations
of non-orientable surfaces were constructed by Ringel~\cite{Ringel1955}.
For minimal triangulations of orientable surfaces see Jungerman and Ringel~\cite{JungermanRingel1980}
and the references contained therein.

\bigskip

A \emph{map} on a surface $M$ is a decomposition of $M$ into\index{map}
polygonal regions or \emph{countries} such that every vertex
has at least degree $3$ and each vertex of degree $s$ is incident
with $s$ countries. If the decomposition is simplicial, the map is called \emph{simplicial}.
For surveys on \emph{polyhedral maps}, i.e., maps that are realized
geometrically in ${\mathbb R}^3$ with straight edges, flat faces,
and without self-intersections, see Brehm and Schulte~\cite{BrehmSchulte1997}
and Brehm and Wills~\cite{BrehmWills1993}.

A map is \emph{regular} if it has a \emph{flag-transitive} automorphism group, 
that is, if its automorphism group is transitive on the triples 
\hbox{$\{\textit{vertex}\subset\textit{edge}\subset\textit{facet}\}$}, or \emph{flags} for short.
Regular maps can be seen as non-spherical combinatorial analogues of the Platonic solids.
Well known examples of flag-transitive simplicial maps with few vertices
are -- besides the tetrahedron, the octahedron, the icosahedron, and ${\mathbb R}{\bf P}^{\,2}_{\,6}$ --
Dyck's regular map \cite{Dyck1880b}, \cite{Dyck1880a} of genus $3$ on $12$ vertices
and Klein's regular map \cite{Klein1879} of genus $3$ on $24$ vertices.
(For polyhedral realizations of Dyck's regular map see Bokowski \cite{Bokowski1989} 
and Brehm~\cite{Brehm1987a}. A polyhedral realization of Klein's map has been 
constructed by Schulte and Wills \cite{SchulteWills1985}.)

A systematic classification of regular maps on orientable surfaces of genus~$1$ and $2$ was
begun by Brahana \cite{Brahana1926}, \cite{Brahana1927} and completed
by Coxeter~\cite{Coxeter1948} and Coxeter and Moser~\cite{CoxeterMoser1957}). 
Coxeter and Moser \cite{CoxeterMoser1957} also constructed regular maps 
on non-orientable surfaces of small genus as well as some interesting
infinite families. Regular maps on the orientable surface of genus $3$ 
were classified by Sherk~\cite{Sherk1959} and 
those on orientable surfaces of genus $4$ to $6$ by Garbe~\cite{Garbe1969}.
Conder and Everitt~\cite{ConderEveritt1995} constructed series 
of regular maps on non-orientable surfaces and Wilson~\cite{Wilson1976}
enumerated all regular maps with up to $100$ edges.
The regular simplicial maps with up to $31$ vertices 
are listed in 
%Chapter~\ref{ch:enumeration}.
\cite{Lutz2004bpre}.

In generalization of regular surfaces, McMullen, Schulz, and Wills \cite{McMullenSchulzWills1982}
called a polyhedral surface \emph{equivelar} of type $\{ p,q\}$ if all its $2$-faces are $p$-gons 
and all its vertices have degree $q$. (In various older papers,
equivelar polyhedra are called regular as well; cf.~\cite{Altshuler1973}.)
Equivelar triangulations of the torus were constructed by Altshuler \cite{Altshuler1973}.
Further series were given by McMullen, Schulz, and Wills \cite{McMullenSchulzWills1982}, \cite{McMullenSchulzWills1983}.
Datta and Nilakantan \cite{DattaNilakantan2001} enumerated all
equivelar surfaces with few vertices, in particular, all simplicial
equivelar polyhedra with up to $11$ vertices. For examples of
face-transitive polyhedra see Wills~\cite{Wills1986}.

\bigskip

\begin{figure}[h]
\begin{center}
\psfrag{1}{\small 1}
\psfrag{2}{\small 2}
\psfrag{3}{\small 3}
\psfrag{4}{\small 4}
\psfrag{5}{\small 5}
\psfrag{6}{\small 6}
\psfrag{7}{\small 7}
\psfrag{Coordinates:}{\small Coordinates:}
\psfrag{1:   (3,-3,0)}{\small 1:   $(3,-3,0)$}
\psfrag{2:   (-3,3,0)}{\small 2:   $(-3,3,0)$}
\psfrag{3:   (-3,-3,1)}{\small 3:   $(-3,-3,1)$}
\psfrag{4:   (3,3,1)}{\small 4:   $(3,3,1)$}
\psfrag{5:   (-1,-2,3)}{\small 5:   $(-1,-2,3)$}
\psfrag{6:   (1,2,3)}{\small 6:   $(1,2,3)$}
\psfrag{7:   (0,0,15)}{\small 7:   $(0,0,15)$}
\psfrag{Triangles:}{\small Triangles:}
\psfrag{123      145      156      345      167      467      247}{\small 123\,\, 145\,\, 156\,\, 345\,\, 167\,\, 467\,\, 247}
\psfrag{124      236      256      346      257      357      137}{\small 124\,\, 236\,\, 256\,\, 346\,\, 257\,\, 357\,\, 137}
\includegraphics[width=.455\linewidth]{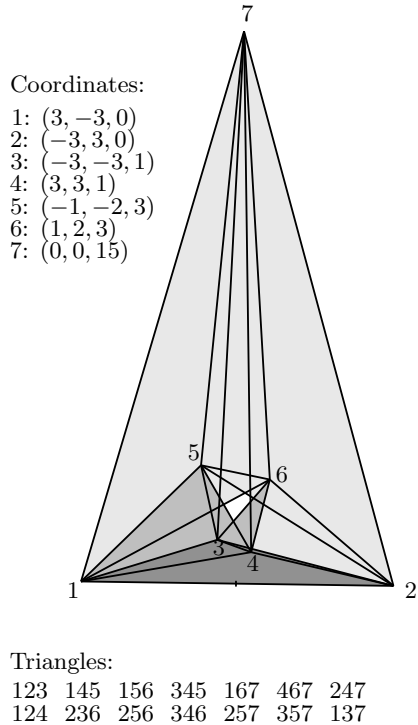}
\end{center}
\caption{{Cs\'asz\'ar}'s torus.}
\label{fig:7torus_small}
\end{figure}

A triangulated surface $M$ with $n$ vertices attains equality in Heawood's bound,
$n=\frac{1}{2}(7+\sqrt{49-24\chi (M)})$,
if and only if the triangulation is \emph{neighborly}, \index{map!neighborly}
that is, the $1$-skeleton of $M$ 
is the complete graph $K_n$. In particular, a triangulation is 
neighborly if and only if $n(7-n)/6=\chi (M)$ is an integer (smaller or
equal to two).
The first examples are the tetrahedron, the real projective plane
${\mathbb R}{\bf P}^{\,2}_6$ with~$6$ vertices, and M\"obius' torus \cite{Moebius1886}
with $7$ vertices.

Neighborly maps with $9$ and $10$ vertices were enumerated by
Altshuler and Brehm \cite{AltshulerBrehm1992}: There are $2$ 
non-orientable neighborly maps of genus $5$ on $9$ vertices 
and $14$ non-orientable neighborly maps of genus $7$ on $10$ vertices. Moreover, there are
$59$ orientable neighborly maps of genus $6$ on $12$ vertices,
as enumerated by Altshuler~\cite{AltshulerBokowskiSchuchert1996} (see also \cite{Altshuler1997}).

Neighborly maps have attracted attention in various ways.
By Theorem~\ref{thm:Heawood}, they provide examples of minimal
triangulations of surfaces. At the same time, every neighborly
map $M$ with $n$ vertices is an example of a minimal graph
embedding of the complete graph $K_n\hookrightarrow M$.

\begin{thm} {\rm (Ringel and Youngs \cite{Ringel1974})}\label{thm:map-thread}
For every surface $M$ different from the $2$-sphere 
and with the exception of the Klein bottle,
there is an embedding\, $K_n\hookrightarrow M$\, if and only if\,
$n\leq\frac{1}{2}(7+\sqrt{49-24\chi (M)})$.
If equality holds, then the embedding of the complete graph induces
a triangulation of $M$. For the Klein bottle, there exist embeddings for\, $n\leq 6$.
\end{thm}

Let the \emph{chromatic number} $\chi_{\rm ch}(M)$ of a surface $M$
be the minimal number of colors that is needed to color any polyhedral
map on $M$. Heawood proved that\, $\chi_{\rm ch}(M)\leq\lfloor \frac{7+\sqrt{1+48g}}{2}\rfloor$\,
for all orientable surfaces of genus $g\geq 1$.

%\pagebreak

As a direct consequence of the existence of embeddings of complete graphs
on surfaces according to Theorem~\ref{thm:map-thread}, Ringel and Youngs \cite{Ringel1974} 
were able to complete in 1968 the \emph{map color theorem},\index{map color theorem}
which first was announced by Heawood~\cite{Heawood1890} in 1890.

\begin{mapthm} {\rm (Ringel and Youngs \cite{Ringel1974})}\label{thm:map-col}
If $M$ is an orientable surface of genus $g\geq 1$, then
\begin{equation}
\chi_{\rm ch}(M)=\left\lfloor \frac{7+\sqrt{1+48g}}{2}\right\rfloor.
\end{equation}
\end{mapthm}
Theorem~\ref{thm:map-thread} implies equality
$\chi_{\rm ch}(M)=\lfloor\frac{1}{2}(7+\sqrt{49-24\chi (M)})\rfloor$
also for all non-orientable surfaces, with the exception of the Klein
bottle $K$, where $\chi_{\rm ch}(K)=6$ \cite{Franklin1934}.
For the $2$-dimensional sphere, $\chi_{\rm ch} (S^2)=4$ is the essence of the famous \emph{Four Color Theorem} 
of Appel and Haken~\cite{AppelHaken1976}; cf.\ also Robertson, Sanders,
Seymour, and Thomas \cite{RobertsonSandersSeymourThomas1997}.

\bigskip

Neighborly maps are of interest also with respect to 
geometric realizations of surfaces.
\begin{ques} {\rm (Gr\"unbaum~\cite[Ch.~13.2]{Gruenbaum1967})}\index{realizability}
Can every triangulated orientable $2$-manifold be embedded 
geometrically in\, ${\mathbb R}^3$, i.e., can it be realized 
with straight edges, flat triangles, and without self intersections?
\end{ques}
By Steinitz' theorem (cf.~\cite{Ziegler1995}), every combinatorial $2$-sphere 
is realizable as the boundary complex of a convex $3$-dimensional polytope. 
For the $2$-torus of genus $1$ the realizability problem is still open.
\begin{conj} {\rm (Duke~\cite{Duke1970})}
Every triangulated torus can be realized as a polyhedron in ${\mathbb R}^3$.
\end{conj}
A first explicit geometric realization of M\"obius' minimal $7$-vertex triangulation\index{torus!Cs\'asz\'ar}
of the $2$-torus was given by Cs\'asz\'ar~\cite{Csaszar1949} (see
Figure~\ref{fig:7torus_small} and \cite{Lutz2002b}).
Bokowski and Eggert \cite{BokowskiEggert1991} showed that
there are altogether $72$ different types of realizations of the M\"obius torus,
and Bokowski and Fendrich \cite{Bokowski-pers} verified that 
triangulated tori with up to $11$ vertices are all realizable.
Sets of coordinates for triangulated tori with up to $10$ vertices
can be found in \cite{HougardyLutzZelke2005bpre}.
%%%%%%%%%%

%%%  -------> n=10

%%%%%%%%%%

Brehm and Bokowski \cite{BokowskiBrehm1987}, \cite{BokowskiBrehm1989}, \cite{Brehm1981}, \cite{Brehm1987b}
constructed geometric realizations for examples of triangulated orientable $2$-manifolds of genus\, $g=2,3,4$\,
with minimal numbers of vertices $n=10,10,11$, respectively.

Neighborly maps of higher genus, however, were considered as candidates for counter-examples 
to the Gr\"unbaum realization problem for a while; cf.~\cite[p.~137]{BokowskiSturmfels1989}. 
Neighborly orientable maps have genus\, $g=(n-3)(n-4)/12$\, and
therefore\, $n\equiv 0,3,4,7\,{\rm mod}\,12$\, vertices, with\, $g=6$\, and\, $n=12$\, as the first
case beyond the tetrahedron and the $7$-vertex torus.

\begin{thm} {\rm (Bokowski and Guedes de Oliveira~\cite{BokowskiGuedes_de_Oliveira2000})}
The triangulated orientable surface $N^{12}_{54}$ of genus $6$ with $12$
vertices of Alshuler's list~\cite{AltshulerBokowskiSchuchert1996}
is not geometrically embeddable in ${\mathbb R}^3$.
\end{thm}

Datta and Nilakantan \cite{DattaNilakantan2002} determined
all triangulated surfaces with $8$~vertices. Those with 
$9$ and $10$ were enumerated in \cite{Lutz2005apre}; 
see Table~\ref{tbl:ten2d_9_10} for the respective numbers
of combinatorial types.

\begin{table}
\small\centering
\defaultaddspace=0.15em
\caption{Triangulated surfaces with up to $10$ vertices.}\label{tbl:ten2d_9_10}
\begin{tabular*}{\linewidth}{@{}l@{\extracolsep{12pt}}l@{\extracolsep{12pt}}r@{\extracolsep{\fill}}l@{\extracolsep{12pt}}l@{\extracolsep{12pt}}r@{\extracolsep{\fill}}l@{\extracolsep{12pt}}l@{\extracolsep{12pt}}r@{}}
 \addlinespace
 \addlinespace
 \addlinespace
 \addlinespace
\toprule
 \addlinespace
 \addlinespace
 \addlinespace
 \addlinespace
  $n$ &   Surface  &  Types  &    $n$ &   Surface  &  Types  &   $n$ &   Surface  &  Types   \\ \cmidrule{1-3}\cmidrule{4-6}\cmidrule{7-9}
 \addlinespace
 \addlinespace
 \addlinespace
 \addlinespace
  4  & $S^2$                  &   1 &   8  & $S^2$                  &  14 &   10  & $S^2$                  &   233 \\
 \addlinespace
     &                        &     &      & $T^2$                  &   7 &       & $T^2$                  &  2109 \\
 \addlinespace
  5  & $S^2$                  &   1 &      &                        &     &       & $M(2,+)$               &   865 \\
 \addlinespace
     &                        &     &      & ${\mathbb R}{\bf P}^2$ &  16 &       & $M(3,+)$               &    20 \\
 \addlinespace
  6  & $S^2$                  &   2 &      & $K^2$                  &   6 &       &                        &       \\
 \addlinespace
     &                        &     &      &                        &     &       & ${\mathbb R}{\bf P}^2$ &  1210 \\
 \addlinespace
     & ${\mathbb R}{\bf P}^2$ &   1 &   9  & $S^2$                  &  50 &       & $K^2$                  &  4462 \\
 \addlinespace
     &                        &     &      & $T^2$                  & 112 &       & $M(3,-)$               & 11784 \\
 \addlinespace
  7  & $S^2$                  &   5 &      &                        &     &       & $M(4,-)$               & 13657 \\
 \addlinespace
     & $T^2$                  &   1 &      & ${\mathbb R}{\bf P}^2$ & 134 &       & $M(5,-)$               &  7050 \\
 \addlinespace
     &                        &     &      & $K^2$                  & 187 &       & $M(6,-)$               &  1022 \\
 \addlinespace
     & ${\mathbb R}{\bf P}^2$ &   3 &      & $M(3,-)$               & 133 &       & $M(7,-)$               &    14 \\
 \addlinespace
     &                        &     &      & $M(4,-)$               &  37 \\
 \addlinespace
     &                        &     &      & $M(5,-)$               &   2 \\
 \addlinespace
 \addlinespace
 \addlinespace
 \addlinespace
\bottomrule
\end{tabular*}
\end{table}

\begin{thm} {\rm (Bokowski and Lutz; cf.~\cite{Lutz2005apre})}
All $865$ vertex-minimal $10$-vertex triangulations of the orientable surface 
of genus $2$ can be realized geometrically in ${\mathbb R}^3$.
%            can be realized as polyhedra in ${\mathbb R}^3$.
\end{thm}

Vertex-transitive triangulations of surfaces with 
up to $15$ vertices are given in
%Chapter~\ref{ch:manifold}.
\cite{KoehlerLutz2004bpre}. 
Enumeration results for vertex-transitive
neighborly triangulations  with up to $22$ vertices
can be found in 
%Chapter~\ref{ch:enumeration}.
\cite{Lutz2004bpre}.

\bigskip

For a given number $n$ of vertices it is way easier to enumerate triangulations of the
$2$-dimensional sphere only than to enumerate all triangulated
$2$-manifolds: According to
Steinitz' theorem (\cite{Steinitz1922}, \cite{SteinitzRademacher1934};
cf.\ \cite{Ziegler1995}), every triangulation of $S^2$ is polytopal
and therefore is, by stereographic projection, equivalent to a planar triangulation (with straight edges).
Br\"uckner \cite{Brueckner1897}, \cite{Brueckner1931} listed (by hand!)
all triangulations of $S^2$ with up to $12$ vertices. His census
was later corrected slightly for $n=11$ by Grace \cite{Grace1965}
and for $n=12$ by Bowen and Fisk \cite{BowenFisk1967}.
Triangulations with up to $23$ vertices were enumerated with the program
\texttt{plantri} by Brinkmann and McKay~\cite{plantri} (see the manual of \texttt{plantri}
and also Royle~\cite{Royle_url}). Table~\ref{tbl:ten2d_spheres} gives the respective numbers.
Precise formulas for rooted triangulations of the $2$-sphere were
determined by Tutte~\cite{Tutte1962}.

\begin{table}
\small\centering
\defaultaddspace=0.15em
\caption{Triangulated $2$-spheres with $11\leq n\leq 23$ vertices.}\label{tbl:ten2d_spheres}
\begin{tabular}{r@{\hspace{20mm}}r}
 \addlinespace
 \addlinespace
 \addlinespace
 \addlinespace
\toprule
 \addlinespace
 \addlinespace
 \addlinespace
  $n$ &   Types \\ \midrule
 \addlinespace
 \addlinespace
 \addlinespace
 \addlinespace
  11  &  1249 \\
 \addlinespace
  12  &  7595 \\
 \addlinespace
  13  &  49566 \\
 \addlinespace
  14  &  339722 \\
 \addlinespace
  15  &  2406841 \\
 \addlinespace
  16  &  17490241 \\
 \addlinespace
  17  &  129664753 \\
 \addlinespace
  18  &  977526957 \\
 \addlinespace
  19  &  7475907149 \\
 \addlinespace
  20  &  57896349553 \\
 \addlinespace
  21  &  453382272049 \\
 \addlinespace
  22  &  3585853662949 \\
 \addlinespace
  23  &  28615703421545 \\
 \addlinespace
 \addlinespace
 \addlinespace
\bottomrule
\end{tabular}
\end{table}

\section{3-Manifolds}
\label{sec:3-manifolds}

All triangulated $3$-manifold with $f$-vector $(n,f_1,f_2,f_3)$ 
satisfy the relations $n - f_1 + f_2 - f_3 = 0$ (Euler) and $2 f_2 = 4
f_3$ (double counting).
Thus
\begin{equation}
f=(n,f_1,2f_1-2n,f_1-n).
\end{equation}
A complete characterization of the $f$-vectors of
the $3$-manifolds $S^3$, $S^2\hbox{$\times\hspace{-1.62ex}\_\hspace{-.4ex}\_\hspace{.7ex}$}S^1$,
$S^2\!\times\!S^1$, and ${\mathbb R}{\bf P}^{\,3}$
was given by Walkup.

\enlargethispage*{3mm}

\begin{thm} {\rm (Walkup \cite{Walkup1970})}\label{thm:Walkup_gamma}
For every $3$-manifold $M$ there is an integer $\gamma (M)$
such that
\begin{equation}
f_1\geq 4n+\gamma (M)
\end{equation}
for every triangulation of $M$ with $n$ vertices and $f_1$ edges. 
Moreover there is an integer $\gamma^*(M)\geq \gamma (M)$
such that for every pair $(n,f_1)$ with $n\geq 0$ and
\begin{equation}
\binom{n}{2}\geq f_1\geq 4n+\gamma^*(M)
\end{equation}
there is a triangulation of $M$ with $n$ vertices and $f_1$ edges.
In particular,
\begin{itemize}
\item[{\rm (a)}] $\gamma^*=\gamma =-10$\, for\, $S^3$,
\item[{\rm (b)}] $\gamma^*=\gamma =0$\, for\, $S^2\hbox{$\times\hspace{-1.62ex}\_\hspace{-.4ex}\_\hspace{.7ex}$}S^1$,
\item[{\rm (c)}] $\gamma^*=1$\, and\, $\gamma =0$\, for\,
  $S^2\!\times\!S^1$, where, with the exception $(9,36)$, all pairs $(n,f_1)$
  with $n\geq 0$ and $4n+\gamma (M)\leq f_1\leq\binom{n}{2}$ occur,
\item[{\rm (d)}] $\gamma^*=\gamma =7$\, for\, ${\mathbb R}{\bf P}^{\,3}$, and
\item[{\rm (e)}] $\gamma^*(M)\geq\gamma (M)\geq 8$\, for all
  other $3$-manifolds $M$.
\end{itemize}
\end{thm}
For an alternative proof of the existence of neighborly triangulations
for every $3$-manifold see Sarkaria~\cite{Sarkaria1983}.

Let us remark that if a $3$-manifold $M$ can be triangulated with $n$ vertices and $f_1$
edges, then triangulations of $M$ with $n+k$ vertices and $f_1+4k$ edges 
can be obtained for $k\geq 0$ by successive \emph{stacking}.
(In every stacking step some tetrahedron of the respective triangulation of $M$ 
is subdivided: This adds one vertex and four edges each.)

\begin{conj}\label{conj:T3}
The $f$-vectors of triangulations of the $3$-torus $T^3$ are characterized
by\, $\gamma^*(T^3)=\gamma (T^3)=45$.
\end{conj}
Conjecture~\ref{conj:T3} implies Conjecture~\ref{conj:Td} for $d=3$.
We used the bistellar flip program BISTELLAR~\cite{Lutz_BISTELLAR} to
verify that there are triangulations of ${\bf T}^3$ for all pairs $(n,f_1)$ 
with\, $15\leq n\leq 35$\, and\, $4n+45\leq f_1\leq\binom{n}{2}$.

Also that there are triangulations of the lens space $L(3,1)$ for all $(n,f_1)$
with\, $12\leq n\leq 35$\, and\, $4n+18\leq f_1\leq\binom{n}{2}$.
\begin{conj}
The $f$-vectors of triangulations of the lens space $L(3,1)$ are characterized
by\, $\gamma^*(L(3,1))=\gamma (L(3,1))=18$.
\end{conj}

%%%%%%%%%%%%%%%%%%

In 1990, K\"uhnel \cite{Kuehnel1990-few} gave a list of six 
(pairwise different) $3$-manifolds for which he knew triangulations
with $15$ or less vertices: The $3$-sphere $S^3$, the twisted $S^2$-bundle over $S^1$ 
(i.e., the $3$-dimensional Klein bottle, which we usually denote by 
$S^2\hbox{$\times\hspace{-1.62ex}\_\hspace{-.4ex}\_\hspace{.7ex}$}S^1$),
the product $S^2\!\times\!S^1$, 
the projective $3$-space ${\mathbb R}{\bf P}^{\,3}$, 
the $3$-dimensional torus ${\bf T}^3$ \cite{KuehnelLassmann1984-3torus}, 
and Cartan's hypersurface $S^3/Q$  (\cite{BrehmKuehnel1986}; \cite{Cartan1939}).\index{Cartan's hypersurface}

Since then, twenty-one further examples have been found.
Brehm and \'Swiatkowski~\cite{BrehmSwiatkowski1993} constructed triangulations of
$L(3,1)$ and $L(4,1)$ with $13$ and $15$ vertices, respectively;
and by local modifications, Brehm~\cite{Brehm-pers} found a triangulation 
of $L(3,1)$ with $12$ vertices. 
K\"uhnel and Lassmann \cite{KuehnelLassmann1985-di} %\linebreak
listed  
two combinatorial $3$-manifolds with $12$ vertices that have homology
$H_*=({\mathbb Z},{\mathbb Z}^2,{\mathbb Z}^2,{\mathbb Z})$
and which are triangulations of $(S^2\!\times\!S^1)\# (S^2\!\times\!S^1)$.
All other eighteen examples are new. In addition, we improved the 
number of vertices for $L(4,1)$ to\, $n=14$.
\begin{thm}\label{thm:less15}
There are at least\, $27$ distinct $3$-manifolds
that can be triangulated with\, $n\leq 15$ vertices; see Table~\ref{tbl:3-list}.
\end{thm}
Triangulations of the respective manifolds were constructed as
described in 
%Chapter~\ref{ch:geometric}
\cite{Lutz2003bpre} 
and \cite{BrehmLutz2002pre}.
To these triangulations we applied bistellar flips
to reduce the numbers of vertices and edges.
%%%%<--------------- Lemma on connected sums!!!!!!!!!!!!!!1 
The resulting triangulations can be found online at \cite{Lutz_PAGE}.

\begin{table}
\small\centering
\defaultaddspace=0.2em
\caption{\protect\parbox[t]{9.5cm}{Combinatorial $3$-manifolds with
    $n\leq 15$ vertices and smallest known transitive triangulations with $n_{vt}$
    vertices (minimal if underlined).}}\label{tbl:3-list}
\begin{tabular*}{\linewidth}{@{\extracolsep{\fill}}llrrl@{}}
 \addlinespace
 \addlinespace
 \addlinespace
\toprule
 \addlinespace
 \addlinespace
   Manifold          & Homology         &        $n$       &      $n_{vt}$    &  Reference \\ \midrule
 \addlinespace
 \addlinespace
 \addlinespace
      $S^3$          & $({\mathbb Z},0,0,{\mathbb Z})$ & $\underline{5}$ &  $\underline{5}$ & \\
 \addlinespace
 $S^2\hbox{$\times\hspace{-1.62ex}\_\hspace{-.4ex}\_\hspace{.7ex}$}S^1$
                     & $({\mathbb Z},{\mathbb Z},{\mathbb Z}_2,0)$ & $\underline{9}$ &  $\underline{9}$ & \cite{AltshulerSteinberg1974}, \cite{Walkup1970} \\
 \addlinespace
 $S^2\!\times\!S^1$  & $({\mathbb Z},{\mathbb Z},{\mathbb Z},{\mathbb Z})$ & $\underline{10}$ & $\underline{10}$ & \cite{Walkup1970} \\
 \addlinespace
 ${\mathbb R}{\bf P}^{\,3}=L(2,1)$
                     & $({\mathbb Z},{\mathbb Z}_2,0,{\mathbb Z})$ & $\underline{11}$ & $\underline{12}$ & \cite{BrehmSwiatkowski1993}, \cite{Walkup1970} \\
 \addlinespace
 $L(3,1)$            & $({\mathbb Z},{\mathbb Z}_3,0,{\mathbb Z})$ & $\underline{12}$ & $\underline{14}$ & \cite{BagchiDatta2003pre}, \cite{Brehm-pers}; \cite{KuehnelLassmann1985-di} \\
 \addlinespace
 $(S^2\!\times\!S^1)\# (S^2\!\times\!S^1)$
                     & $({\mathbb Z},{\mathbb Z}^2,{\mathbb Z}^2,{\mathbb Z})$ &               12 & $\underline{12}$ & \cite{KuehnelLassmann1985-di} \\ 
 $=(S^2\!\times\!S^1)\# -(S^2\!\times\!S^1)$ &&&& \\
 \addlinespace
 $(S^2\hbox{$\times\hspace{-1.62ex}\_\hspace{-.4ex}\_\hspace{.7ex}$}S^1)\# (S^2\hbox{$\times\hspace{-1.62ex}\_\hspace{-.4ex}\_\hspace{.7ex}$}S^1)$
                     & $({\mathbb Z},{\mathbb Z}^2,{\mathbb Z}\oplus{\mathbb Z}_2,0)$ &               12 & $\underline{16}$ & \\
 $=(S^2\hbox{$\times\hspace{-1.62ex}\_\hspace{-.4ex}\_\hspace{.7ex}$}S^1)\# (S^2\!\times\!S^1)$ &&&& \\
 \addlinespace
 $(S^2\!\times\!S^1)^{\# 3}$ & $({\mathbb Z},{\mathbb Z}^3,{\mathbb Z}^3,{\mathbb Z})$ &               13 & $\underline{16}$ & \\
 \addlinespace
 $(S^2\hbox{$\times\hspace{-1.62ex}\_\hspace{-.4ex}\_\hspace{.7ex}$}S^1)^{\# 3}$
                     & $({\mathbb Z},{\mathbb Z}^3,{\mathbb Z}^2\oplus{\mathbb Z}_2,0)$ &               13 & $\underline{16}$ & \\
 \addlinespace
 $L(4,1)$            & $({\mathbb Z},{\mathbb Z}_4,0,{\mathbb Z})$ &               14 & $\underline{16}$ & \\
 \addlinespace
 $L(5,2)$            & $({\mathbb Z},{\mathbb Z}_5,0,{\mathbb Z})$ &               14 &             ? & \\
 \addlinespace
 $(S^2\!\times\!S^1)\#\,{\mathbb R}{\bf P}^{\,3}$ & $({\mathbb Z},{\mathbb Z}\oplus{\mathbb Z}_2,{\mathbb Z},{\mathbb Z})$ &               14 &     ? & \\
 $=(S^2\!\times\!S^1)\# -{\mathbb R}{\bf P}^{\,3}$ &&&& \\
 \addlinespace
 $(S^2\hbox{$\times\hspace{-1.62ex}\_\hspace{-.4ex}\_\hspace{.7ex}$}S^1)\#\,{\mathbb R}{\bf P}^{\,3}$
                     & $({\mathbb Z},{\mathbb Z}\oplus{\mathbb Z}_2,{\mathbb Z}_2,0)$ &               14 &            ? & \\
 \addlinespace
 ${\mathbb R}{\bf P}^{\,2}\!\times\!S^1$
                     & $({\mathbb Z},{\mathbb Z}\oplus{\mathbb Z}_2,{\mathbb Z}_2,0)$ &               14 &               $\underline{17}$ & \cite{KuehnelLassmann1985-di} \\
 \addlinespace
 $(S^2\hbox{$\times\hspace{-1.62ex}\_\hspace{-.4ex}\_\hspace{.7ex}$}S^1)^{\# 2}\#\,{\mathbb R}{\bf P}^{\,3}$
                     & $({\mathbb Z},{\mathbb Z}^2\oplus{\mathbb Z}_2,{\mathbb Z}\oplus{\mathbb Z}_2,0)$ &               14 &            ? & \\
 \addlinespace
 $S^3/T^*$           & $({\mathbb Z},{\mathbb Z}_3,0,{\mathbb Z})$ &               15 & $\underline{16}$ & \\
 \addlinespace
 $L(5,1)$            & $({\mathbb Z},{\mathbb Z}_5,0,{\mathbb Z})$ &               15 &            ? & \\
 \addlinespace
 $P(2)=S^3/Q$        & $({\mathbb Z},{\mathbb Z}_2^{\,2},0,{\mathbb Z})$ &         15 & $\underline{15}$ & \cite{BrehmKuehnel1986} \\
 \addlinespace
 $P(3)$              & $({\mathbb Z},{\mathbb Z}_4,0,{\mathbb Z})$ &               15 &            ? & \\
 \addlinespace
 $P(4)$              & $({\mathbb Z},{\mathbb Z}_2^{\,2},0,{\mathbb Z})$ &         15 &            ? & \\
 \addlinespace
 ${\mathbb R}{\bf P}^{\,3}\#\,{\mathbb R}{\bf P}^{\,3}$ & $({\mathbb Z},{\mathbb Z}_2^{\,2},0,{\mathbb Z})$ &               15 &     ? & \\
 $={\mathbb R}{\bf P}^{\,3}\# -{\mathbb R}{\bf P}^{\,3}$ &&&& \\
 \addlinespace
 $(S^2\!\times\!S^1)\# L(3,1)$
                     & $({\mathbb Z},{\mathbb Z}\oplus{\mathbb Z}_3,{\mathbb Z},{\mathbb Z})$ &       15 &            ? & \\
 \addlinespace
 $(S^2\hbox{$\times\hspace{-1.62ex}\_\hspace{-.4ex}\_\hspace{.7ex}$}S^1)\# L(3,1)$
                     & $({\mathbb Z},{\mathbb Z}\oplus{\mathbb Z}_3,{\mathbb Z}_2,0)$ &               15 &            ? & \\
 \addlinespace
 $(S^2\!\times\!S^1)^{\# 2}\#\,{\mathbb R}{\bf P}^{\,3}$
                     & $({\mathbb Z},{\mathbb Z}^2\oplus{\mathbb Z}_2,{\mathbb Z}^2,{\mathbb Z})$ &   15 &            ? & \\
 \addlinespace
 ${\bf T}^3=S^1\!\times\!S^1\!\times\!S^1$ & $({\mathbb Z},{\mathbb Z}^3,{\mathbb Z}^3,{\mathbb Z})$ &               15 & $\underline{15}$ & \cite{KuehnelLassmann1984-3torus} \\
 \addlinespace
 $(S^2\!\times\!S^1)^{\# 4}$ & $({\mathbb Z},{\mathbb Z}^4,{\mathbb Z}^4,{\mathbb Z})$ &              15 &            ? & \\
 \addlinespace
 $(S^2\hbox{$\times\hspace{-1.62ex}\_\hspace{-.4ex}\_\hspace{.7ex}$}S^1)^{\# 4}$
                     & $({\mathbb Z},{\mathbb Z}^4,{\mathbb Z}^3\oplus{\mathbb Z}_2,0)$ &             15 &            ? & \\
% \addlinespace
% $\left. \begin{array}{r}
%         (S^2\!\times\!S^1)\# L(3,1)\\
%         (S^2\!\times\!S^1)\# -L(3,1)
%         \end{array}\right\}$ & $({\mathbb Z},{\mathbb Z}\oplus{\mathbb Z}_3,{\mathbb Z},{\mathbb Z})$ &               15 &            ? & \\ [2.8mm]\hline
 \addlinespace
 \addlinespace
\bottomrule
\end{tabular*}
\end{table}

\begin{conj}\label{conj:less14}
There are only nine $3$-manifolds that can minimally be triangulated
with $n\leq 13$ vertices: $S^3$ with $5$, $S^2\hbox{$\times\hspace{-1.62ex}\_\hspace{-.4ex}\_\hspace{.7ex}$}S^1$
with $9$, $S^2\!\times\!S^1$ with $10$, ${\mathbb R}{\bf P}^{\,3}$
with $11$, $L(3,1)$, $(S^2\!\times\!S^1)^{\#2}$, and $(S^2\hbox{$\times\hspace{-1.62ex}\_\hspace{-.4ex}\_\hspace{.7ex}$}S^1)^{\#2}$
with $12$, as well as $(S^2\!\times\!S^1)^{\#3}$ and $(S^2\hbox{$\times\hspace{-1.62ex}\_\hspace{-.4ex}\_\hspace{.7ex}$}S^1)^{\#3}$
with $13$ vertices. For all other $3$-manifolds at least $14$ vertices are needed.
\end{conj}
In view of Walkup's Theorem~\ref{thm:Walkup_gamma} and Bagchi and Datta's Theorem~\ref{thm:homology_sphere},
the open part of the conjecture (see also Conjecture~\ref{conj:mintri}(a) and (b)) 
is to show that $(S^2\!\times\!S^1)^{\#3}$ and $(S^2\hbox{$\times\hspace{-1.62ex}\_\hspace{-.4ex}\_\hspace{.7ex}$}S^1)^{\#3}$
cannot be triangulated with $12$ vertices and that
$(S^2\!\times\!S^1)^{\#2}$ and $(S^2\hbox{$\times\hspace{-1.62ex}\_\hspace{-.4ex}\_\hspace{.7ex}$}S^1)^{\#2}$
cannot be triangulated with $11$ vertices, while all other
$3$-manifolds, different from the nine listed ones, need at least $14$
vertices for a triangulation. K\"uhnel's Conjecture~\ref{conj:kuehnel_pascal}
would imply that $13$ vertices is best possible for 
$(S^2\!\times\!S^1)^{\#3}$ and $(S^2\hbox{$\times\hspace{-1.62ex}\_\hspace{-.4ex}\_\hspace{.7ex}$}S^1)^{\#3}$.

From Conjecture~\ref{conj:less14} it would follow 
that the $14$-vertex triangulations of $L(4,1)$, $L(5,1)$, $(S^2\!\times\!S^1)\#\,{\mathbb R}{\bf P}^{\,3}$,
$(S^2\hbox{$\times\hspace{-1.62ex}\_\hspace{-.4ex}\_\hspace{.7ex}$}S^1)\#\,{\mathbb R}{\bf P}^{\,3}$,
${\mathbb R}{\bf P}^{\,2}\!\times\!S^1$, and of the connected sum
$(S^2\hbox{$\times\hspace{-1.62ex}\_\hspace{-.4ex}\_\hspace{.7ex}$}S^1)^{\# 2}\#\,{\mathbb R}{\bf P}^{\,3}$
from Theorem~\ref{thm:less15} are vertex-minimal (which comprises Conjecture~\ref{conj:mintri}(a)).

\bigskip

Among the triangulations from Table~\ref{tbl:3-list} are various
examples that have the same homology. We list these examples together
with their fundamental groups and their smallest known $f$-vectors 
in Table~\ref{tbl:hom-iso}.

\begin{table}
\small\centering
\defaultaddspace=0.2em
\caption{$3$-manifolds with isomorphic homology groups on $n\leq 15$ vertices.}\label{tbl:hom-iso}
\begin{tabular*}{\linewidth}{@{\extracolsep{\fill}}llll@{}}
 \addlinespace
 \addlinespace
 \addlinespace
\toprule
 \addlinespace
 \addlinespace
    Manifold         & Homology         & $\pi_1$ &   $f$-vector   \\ \midrule
 \addlinespace
 \addlinespace
$L(3,1)$     & $({\mathbb Z},{\mathbb Z}_3,0,{\mathbb Z})$ & ${\mathbb Z}_3$ & (12,\underline{66},108,54) \\
 \addlinespace
$S^3/T^*$    & $({\mathbb Z},{\mathbb Z}_3,0,{\mathbb Z})$  & $T^*$  & (15,102,174,87) \\ \midrule
 \addlinespace
 \addlinespace
$(S^2\!\times\!S^1)^{\# 3}$ & $({\mathbb Z},{\mathbb Z}^3,{\mathbb Z}^3,{\mathbb Z})$ & ${\mathbb Z}*{\mathbb Z}*{\mathbb Z}$ & (13,72,118,59) \\
 \addlinespace
${\bf T}^3$  & $({\mathbb Z},{\mathbb Z}^3,{\mathbb Z}^3,{\mathbb Z})$  & ${\mathbb Z}^3$  & (15,\underline{105},180,90) \\ \midrule
 \addlinespace
 \addlinespace
$(S^2\hbox{$\times\hspace{-1.62ex}\_\hspace{-.4ex}\_\hspace{.7ex}$}S^1)\#\,{\mathbb R}{\bf P}^{\,3}$ & $({\mathbb Z},{\mathbb Z}\oplus{\mathbb Z}_2,{\mathbb Z}_2,0)$ & ${\mathbb Z}*{\mathbb Z}_2$ & (14,73,118,59) \\
 \addlinespace
${\mathbb R}{\bf P}^{\,2}\!\times\!S^1$ & $({\mathbb Z},{\mathbb Z}\oplus{\mathbb Z}_2,{\mathbb Z}_2,0)$ & ${\mathbb Z}\oplus{\mathbb Z}_2$ & (14,84,140,70) \\ \midrule
 \addlinespace
 \addlinespace
$L(4,1)$ & $({\mathbb Z},{\mathbb Z}_4,0,{\mathbb Z})$ & ${\mathbb Z}_4$  & (14,84,140,70) \\
 \addlinespace
$P(3)$ & $({\mathbb Z},{\mathbb Z}_4,0,{\mathbb Z})$ & $D^*_{3}$  & (15,97,164,82) \\\midrule
 \addlinespace
 \addlinespace
$L(5,2)$ & $({\mathbb Z},{\mathbb Z}_5,0,{\mathbb Z})$ & ${\mathbb Z}_5$  & (14,87,146,73) \\
 \addlinespace
$L(5,1)$ & $({\mathbb Z},{\mathbb Z}_5,0,{\mathbb Z})$ & ${\mathbb Z}_5$  & (15,97,164,82) \\\midrule
 \addlinespace
 \addlinespace
${\mathbb R}{\bf P}^{\,3}\#\,{\mathbb R}{\bf P}^{\,3}$ & $({\mathbb Z},{\mathbb Z}_2^{\,2},0,{\mathbb Z})$ & ${\mathbb Z}_2*{\mathbb Z}_2$ & (15,86,142,71) \\
 \addlinespace
$P(2)=S^3/Q$ & $({\mathbb Z},{\mathbb Z}_2^{\,2},0,{\mathbb Z})$ & $D^*_{2}=Q$  & (15,90,150,75) \\
 \addlinespace
$P(4)$ & $({\mathbb Z},{\mathbb Z}_2^{\,2},0,{\mathbb Z})$ & $D^*_{4}$  & (15,104,178,89) \\
 \addlinespace
\bottomrule
\end{tabular*}
\end{table}

Suppose, one of the manifolds from Table~\ref{tbl:hom-iso} is given to us
as a simplicial complex without further information.
Then computing its homology vector, fundamental group, 
and (as a `quasi-invariant') the $f$-vector of the smallest triangulation that we achieve
from the given complex by bistellar flips, will allow us to make
a quite accurate guess for its topological type.
In general, however, many manifolds will share the same minimal $f$-vector.

The smallest $f$-vectors that we found for $k$-fold connected sums of sphere products $(S^2\!\times\!S^1)^{\# k}$
and twisted sphere products $(S^2\hbox{$\times\hspace{-1.62ex}\_\hspace{-.4ex}\_\hspace{.7ex}$}S^1)^{\# k}$
are identical for $2\leq k\leq 5$; see Table~\ref{tbl:small-connected-sums}.

\begin{table}
\small\centering
\defaultaddspace=0.2em
\caption{Smallest known triangulations of connected sums of $S^2\!\times\!S^1$ and $S^2\hbox{$\times\hspace{-1.62ex}\_\hspace{-.4ex}\_\hspace{.7ex}$}S^1$.}\label{tbl:small-connected-sums}
\begin{tabular}{@{}l@{\hspace{15mm}}l@{}}
 \addlinespace
 \addlinespace
 \addlinespace
\toprule
 \addlinespace
 \addlinespace
    Manifold         & $f$-Vector \\ \midrule
 \addlinespace
 \addlinespace
 \addlinespace
 $(S^2\!\times\!S^1)^{\# 2}$,\, $(S^2\hbox{$\times\hspace{-1.62ex}\_\hspace{-.4ex}\_\hspace{.7ex}$}S^1)^{\# 2}$  & (12,58,92,46) \\
 \addlinespace
 $(S^2\!\times\!S^1)^{\# 3}$,\, $(S^2\hbox{$\times\hspace{-1.62ex}\_\hspace{-.4ex}\_\hspace{.7ex}$}S^1)^{\# 3}$  & (13,72,118,59) \\
 \addlinespace
 $(S^2\!\times\!S^1)^{\# 4}$,\, $(S^2\hbox{$\times\hspace{-1.62ex}\_\hspace{-.4ex}\_\hspace{.7ex}$}S^1)^{\# 4}$  & (15,90,150,75) \\
 \addlinespace
 $(S^2\!\times\!S^1)^{\# 5}$,\, $(S^2\hbox{$\times\hspace{-1.62ex}\_\hspace{-.4ex}\_\hspace{.7ex}$}S^1)^{\# 5}$  & (16,104,176,88) \\
 \addlinespace
 \addlinespace
\bottomrule
\end{tabular}
\end{table}

\begin{thm}  {\rm (Brehm and \'Swiatkowski \cite{BrehmSwiatkowski1993})}
The number of topologically distinct lens spaces
that can be triangulated with $n$ vertices grows
exponentially with $n$.
\end{thm}
Brehm and \'Swiatkowski \cite{BrehmSwiatkowski1993} constructed explicit
triangulations of all lens spaces $L(p,q)$. In particular, they gave an
infinite series of $D_{2(p+2)}$-symmetric triangulations
$S_{2(p+2)}$ of\, $L(p,1)$\, with\, $2p+7$\, vertices.

The Brehm and \'Swiatkowski example $S_{2\cdot 4}$ with
$D_{8}$-symmetry on $11$-vertices is combinatorially 
isomorphic to Walkup's minimal triangulation \cite{Walkup1970} 
of the lens space ${\mathbb R}{\bf P}^3=L(2,1)$.
The facets of this triangulation ${\mathbb R}{\bf P}^3_{11}$ 
with $f$-vector $(11,51,80,40)$ are:

{\small
\begin{center}
\begin{tabular}{llllllll}
$1237$    & $123\,11$ & $1269$    & $126\,11$ & $1279$    & $135\,10$ & $135\,11$ & $137\,10$ \\
$1479$    & $147\,10$ & $1489$    & $148\,10$ & $1568$    & $156\,11$ & $158\,10$ & $1689$ \\
$2348$    & $234\,11$ & $2378$    & $246\,10$ & $246\,11$ & $248\,10$ & $2578$    & $2579$ \\
$258\,10$ & $259\,10$ & $269\,10$ & $3459$    & $345\,11$ & $3489$    & $359\,10$ & $3678$ \\
$367\,10$ & $3689$    & $369\,10$ & $4567$    & $456\,11$ & $4579$    & $467\,10$ & $5678$.
\end{tabular}
\end{center}
}

\noindent
The full automorphism group of ${\mathbb R}{\bf P}^3_{11}$\index{projective space!real}
is larger than $D_{8}$: For every vertex, we computed the
Altshuler-Steinberg determinant \cite{AltshulerSteinberg1973}\,
$\det AA^T$ of the vertex-facet incidence matrix $A$
of the respective vertex-link. Vertices 1--6 yield Altshuler-Steinberg determinant $41616$, 
the determinant for vertices 7--10 is $12096$, and vertex $11$ gives determinant $0$. 
Thus, the automorphism group of ${\mathbb R}{\bf P}^3_{11}$ must be a subgroup of $S_6\times S_4$. 
In fact, it is $2S_4$ with generators\, (1,2,3,4,5,6)(7,8,9),\, (1,2)(3,6)(4,5)(7,9),\, and\, (3,6)(7,9)(8,10), 
which can easily be verified by a computer.

By applying bistellar flips to the triangulation $S_{2\cdot 5}$
with $13$ vertices, we obtained a $12$-vertex triangulations
$L(3,1)_{12}$ of the lens space $L(3,1)$ with $f$-vector 
$(12,\underline{66},108,54)$ and facets

{\small
\begin{center}
\begin{tabular}{lllllll}
$1234$       & $123\,10$    & $1249$       & $1256$       & $1259$       & $126\,11$    & $12\,10\,11$ \\
$1347$       & $1378$       & $138\,10$    & $1479$       & $156\,12$    & $1579$       & $157\,12$    \\
$16\,11\,12$ & $178\,12$    & $18\,10\,11$ & $18\,11\,12$ & $234\,12$    & $23\,10\,12$ & $2489$       \\
$248\,12$    & $2568$       & $2589$       & $2678$       & $267\,11$    & $278\,12$    & $27\,10\,11$ \\
$27\,10\,12$ & $3456$       & $345\,11$    & $3467$       & $34\,11\,12$ & $3568$       & $3589$       \\
$359\,11$    & $3678$       & $389\,10$    & $39\,10\,12$ & $39\,11\,12$ & $456\,10$    & $45\,10\,11$ \\
$4679$       & $469\,10$    & $489\,10$    & $48\,10\,11$ & $48\,11\,12$ & $56\,10\,12$ & $579\,11$    \\
$57\,10\,11$ & $57\,10\,12$ & $679\,11$    & $69\,10\,12$ & $69\,11\,12$.
\end{tabular}
\end{center}
}

\noindent
The symmetry group of $L(3,1)_{12}$ is $S_3$ as a subgroup of $S_6\times S_3\times S_3$
with generators  (1,2)(3,6)(4,5)(7,8)(10,11) and (1,3,5)(2,4,6)(7,8,9)(10,11,12).
(The Altshuler-Steinberg determinant is 134784 for the vertices 1--6,
133056 for the vertices 7--9, and 112320 for the vertices 10--12.)

As already mentioned, Brehm \cite{Brehm-pers} previously 
found a triangulation of $L(3,1)$ with $12$ vertices 
by modifying $S_{2\cdot 5}$. 

\begin{conj}
The minimal triangulation ${\mathbb R}{\bf P}^3_{11}$ of\, ${\mathbb R}{\bf P}^3$ 
is unique with $f$-vector $(11,51,80,40)$.
Also, the minimal triangulation $L(3,1)_{12}$ of $L(3,1)$
is unique with $f$-vector $(12,\underline{66},108,54)$
and is the only triangulation of $L(3,1)$
with $12$ vertices.
\end{conj}

%\pagebreak

The exact numbers of different combinatorial types of triangulations
of $S^3$, $S^2\hbox{$\times\hspace{-1.62ex}\_\hspace{-.4ex}\_\hspace{.7ex}$}S^1$,
and $S^2\!\times\!S^1$ with up to $10$ vertices were obtained by

\begin{tabbing}
      ooooooooooooooooppppppppp\hspace{1mm} \=  \kill
      Gr\"unbaum and Sreedharan \cite{GruenbaumSreedharan1967} \> (simplicial $4$-polytopes with $8$ vertices),\\
      Barnette \cite{Barnette1973c}                            \> (combinatorial $3$-spheres with $8$ vertices),\\
      Altshuler and Steinberg \cite{AltshulerSteinberg1973}    \> (neighborly $4$-polytopes with $9$ vertices),\\
      Altshuler and Steinberg \cite{AltshulerSteinberg1974}    \> (neighborly $3$-manifolds with $9$ vertices),\\
      Altshuler and Steinberg \cite{AltshulerSteinberg1976}    \> (combinatorial $3$-manifolds with $9$ vertices),\\
      Altshuler \cite{Altshuler1977}                           \> (neighborly $3$-manifolds with $10$ vertices), \\
      Lutz \cite{BokowskiBremnerLutzMartin2003pre}             \> (combinatorial $3$-manifolds with $10$ vertices).
\end{tabbing}
For a discussion of the polytopality of the simplicial $3$-spheres with $9$ 
vertices see Altshuler, Bokowski, and Steinberg \cite{AltshulerBokowskiSteinberg1980}
and Engel \cite{Engel1991}. For the polytopality of the neighborly simplicial $3$-spheres with $10$ 
vertices see Bokowski and Garms \cite{BokowskiGarms1987}
and Bokowski and Sturmfels \cite{BokowskiSturmfels1987}. 
The numbers of combinatorial types of $3$-manifolds with up to $10$
vertices can be found in Table~\ref{tbl:ten3d_9_10}.

\begin{table}
\small\centering
\defaultaddspace=0.15em
\caption{Triangulated $3$-manifolds with up to $10$ vertices.}\label{tbl:ten3d_9_10}
\begin{tabular}{r@{\hspace{10mm}}l@{\hspace{10mm}}r@{\hspace{10mm}}r}
 \addlinespace
 \addlinespace
 \addlinespace
 \addlinespace
\toprule
 \addlinespace
 \addlinespace
 \addlinespace
  $n$ &  Manifold          &  Types (all) &   Types (neighborly)\\ \midrule
 \addlinespace
 \addlinespace
 \addlinespace
 \addlinespace
  5  & $S^3$                  & 1    & 1 \\
 \addlinespace
  6  & $S^3$                  & 2    & 1 \\
 \addlinespace
  7  & $S^3$                  & 5    & 1 \\
 \addlinespace
  8  & $S^3$                  & 39   & 4 \\
 \addlinespace
  9  & $S^3$                  & 1296 & 50 \\
 \addlinespace
     & $S^2\hbox{$\times\hspace{-1.62ex}\_\hspace{-.4ex}\_\hspace{.7ex}$}S^1$  & 1 & 1 \\
 \addlinespace
 10  & $S^3$                  & 247882 & 3540 \\
 \addlinespace
     & $S^2\hbox{$\times\hspace{-1.62ex}\_\hspace{-.4ex}\_\hspace{.7ex}$}S^1$  & 615 & 83\\
 \addlinespace
     & $S^2\!\times\!S^1$     & 518 & 54 \\
 \addlinespace
 \addlinespace
 \addlinespace
\bottomrule
\end{tabular}
\end{table}

An upper bound on the numbers of combinatorial types of simplicial 
$4$-polytopes was given by Goodman and Pollack
\cite{GoodmanPollack1986a}, \cite{GoodmanPollack1986b}.
Many combinatorially different types of triangulated spheres 
for growing $n$ were constructed by Kalai~\cite{Kalai1988} 
(for $d\geq 4$) and by Pfeifle and Ziegler \cite{PfeifleZiegler2002pre} (for $d=3$).

K\"uhnel and Lassmann \cite{KuehnelLassmann1985-di} enumerated  all
combinatorial $3$-manifolds with $n\leq 15$ vertices that have a
vertex-transitive cyclic group action as well as all $3$-manifolds 
with a vertex-transitive dihedral action for $n\leq 19$.
(In their list appear two non-orientable manifolds ${\rm IV}_{17}$ and ${\rm IV}_{19}$.
These are homeomorphic to ${\mathbb R}{\bf P}^{\,2}\!\times\!S^1$,
as we were able to recognize with bistellar flips.)
Enumeration results of all vertex-transitive triangulations of $3$-manifolds
with up to $15$ vertices are listed in 
%Chapter~\ref{ch:manifold} 
\cite{KoehlerLutz2004bpre}
and with $16$ and $17$ vertices in 
%Chapter~\ref{ch:enumeration}.
\cite{Lutz2004bpre}.

%{\sc Remark:} In Chapter~\ref{ch:diagonal}, % <---
%we describe two non-orientable 3-manifolds on 16 vertices with
%homology $({\mathbb Z},{\mathbb Z}^2,{\mathbb Z}\oplus{\mathbb Z}_2,0)$ 
%for which it has not been possible to find a triangulation with fewer vertices. 
%Hence, these manifolds might be distinct from \label{misc:not-sum}
%$(S^2\hbox{$\times\hspace{-1.62ex}\_\hspace{-.4ex}\_\hspace{.7ex}$}S^1)\# (S^2\hbox{$\times\hspace{-1.62ex}\_\hspace{-.4ex}\_\hspace{.7ex}$}S^1)$
%which can be triangulated with 12 vertices.
%\bigskip

%%%%%%%%%%%%%%%%%%

\section{4-Manifolds}
\label{sec:4-manifolds}

The unique $9$-vertex triangulation ${\mathbb C}{\bf P}^{\,2}_9$ of K\"uhnel \cite{KuehnelBanchoff1983}
of the complex projective plane ${\mathbb C}{\bf P}^{\,2}$ certainly is the most prominent 
combinatorial $4$-manifold. \index{projective space!complex}
By the Brehm-K\"uhnel bound (Theorem~\ref{thm:bk-minimal}(a), \cite{BrehmKuehnel1987}), 
it has the minimal number of vertices that a combinatorial
$4$-manifold, different from~$S^4$, can have.

K\"uhnel's triangulation $M^4=M^4_{3}(11)$ (\cite{Kuehnel1986a-series}, \cite{KuehnelLassmann1996-bundle}) 
of the product $S^3\!\times\!S^1$ is vertex-minimal with $11$ vertices.
(The combinatorial manifolds $M^4_{3}(n)$ of K\"uhnel and Lassmann \cite{KuehnelLassmann1996-bundle}
are triangulations of $S^3\!\times\!S^1$ for all $n\geq 11$.)

For the twisted sphere product $S^3\hbox{$\times\hspace{-1.62ex}\_\hspace{-.4ex}\_\hspace{.7ex}$}S^1$
it is conjectured (Conjecture~\ref{conj:products_even_dim}) that at least
$12$ vertices are needed for a triangulation, while $11$ vertices is
the current best lower bound according to Corollary~\ref{cor:BrehmKuehnel}.
A vertex-transitive $12$-vertex triangulation $\mbox{}^4\hspace{.3pt}12^{\,54}_{\,1}$ of
$S^3\hbox{$\times\hspace{-1.62ex}\_\hspace{-.4ex}\_\hspace{.7ex}$}S^1$
with $f$-vector $(12,\underline{66},144,150,60)$ is described in 
%Chapter~\ref{ch:manifold}.
\cite{KoehlerLutz2004bpre}.
We applied the bistellar flip program BISTELLAR~\cite{Lutz_BISTELLAR} to the
product triangulation of $S^3\hbox{$\times\hspace{-1.62ex}\_\hspace{-.4ex}\_\hspace{.7ex}$}S^1$
from Proposition~\ref{prop:series_2d+3}:
\begin{prop}
There is a $12$-vertex triangulation $(S^3\hbox{$\times\hspace{-1.62ex}\_\hspace{-.4ex}\_\hspace{.7ex}$}S^1)_{12}$
of\, $S^3\hbox{$\times\hspace{-1.62ex}\_\hspace{-.4ex}\_\hspace{.7ex}$}S^1$
with $f$-vector $(12,60,120,120,48)$.
\end{prop}
The facets of
$(S^3\hbox{$\times\hspace{-1.62ex}\_\hspace{-.4ex}\_\hspace{.7ex}$}S^1)_{12}$ are
\begin{center}\small
\begin{tabular}{l@{\hspace{4mm}}l@{\hspace{4mm}}l@{\hspace{4mm}}l@{\hspace{4mm}}l@{\hspace{4mm}}l@{\hspace{4mm}}l@{\hspace{4mm}}l}
 $12678$ & $12679$ & $1268a$ & $1269a$ & $1278c$ & $1279a$ & $127ac$ & $128ac$ \\
 $1348b$ & $1348c$ & $134bc$ & $1358b$ & $1358c$ & $135bc$ & $148bc$ & $1578a$ \\
 $1578b$ & $157ac$ & $157bc$ & $158ac$ & $1678a$ & $1679a$ & $178bc$ & $23467$ \\ 
 $2346b$ & $23479$ & $2349b$ & $2367a$ & $2369a$ & $2369b$ & $2379a$ & $24679$ \\
 $2469b$ & $2678a$ & $278ac$ & $34679$ & $3469c$ & $346bc$ & $3489b$ & $3489c$ \\
 $358bc$ & $3679a$ & $369bc$ & $389bc$ & $469bc$ & $489bc$ & $578ac$ & $578bc$,
\end{tabular}
\end{center}
with vertices $1$--$9$, $a$, $b$, and $c$, respectively.
\begin{conj}
The $f$-vector $(12,60,120,120,48)$ is com\-po\-nent-wise minimal
for combinatorial triangulations of\,
$S^3\hbox{$\times\hspace{-1.62ex}\_\hspace{-.4ex}\_\hspace{.7ex}$}S^1$.
\end{conj}

In Table~\ref{tbl:4-list} we list the $4$-dimensional manifolds for
which we know triangulations with $n\leq 16$ vertices.

\begin{table}
\small\centering
\defaultaddspace=0.2em
\caption{\protect\parbox[t]{9.5cm}{Combinatorial $4$-manifolds with
    $n\leq 16$ vertices and smallest known transitive triangulations with $n_{vt}$
    vertices (minimal if underlined).}}\label{tbl:4-list}
\begin{tabular*}{\linewidth}{@{\extracolsep{\fill}}llrrl@{}}
 \addlinespace
 \addlinespace
 \addlinespace
\toprule
 \addlinespace
 \addlinespace
   Manifold          & Homology         &        $n$       &      $n_{vt}$    &  Reference \\ \midrule
 \addlinespace
 \addlinespace
 \addlinespace
     $S^4$           & $({\mathbb Z},0,0,0,{\mathbb Z})$ &    $\underline{6}$ & $\underline{6}$ & \\
 \addlinespace
 ${\mathbb C}{\bf P}^{\,2}$
                     & $({\mathbb Z},0,{\mathbb Z},0,{\mathbb Z})$ &     $\underline{9}$ & $\underline{9}$ & \cite{KuehnelBanchoff1983}, \cite{KuehnelLassmann1983-unique} \\
 \addlinespace
 $S^3\!\times\!S^1$  & $({\mathbb Z},{\mathbb Z},0,{\mathbb Z},{\mathbb Z})$ &    $\underline{11}$ & $\underline{11}$ & \cite{Kuehnel1986a-series} \\
 \addlinespace
 $S^2\!\times\!S^2$  & $({\mathbb Z},0,{\mathbb Z}^2,0,{\mathbb Z})$ &    $\underline{11}$ & $\underline{12}$ & \cite{Sparla1997}, \cite{Sparla1998}\\ %  Ch.~\ref{ch:centrally}; 
 \addlinespace
 $S^3\hbox{$\times\hspace{-1.62ex}\_\hspace{-.4ex}\_\hspace{.7ex}$}S^1$
                     & $({\mathbb Z},{\mathbb Z},0,{\mathbb Z}_2,0)$ &    12 & $\underline{12}$ & \cite{KoehlerLutz2004bpre}\\%Ch.~\ref{ch:manifold}\\
 \addlinespace
 ${\mathbb C}{\bf P}^{\,2}\#\,{\mathbb C}{\bf P}^{\,2}$
                     & $({\mathbb Z},0,{\mathbb Z}^2,0,{\mathbb Z})$ &    12 & ? & \\
 \addlinespace
 ${\mathbb C}{\bf P}^{\,2}\# -{\mathbb C}{\bf P}^{\,2}=S^2\hbox{$\times\hspace{-1.62ex}\_\hspace{-.4ex}\_\hspace{.7ex}$}S^2$
                     & $({\mathbb Z},0,{\mathbb Z}^2,0,{\mathbb Z})$ &    12 & ? & \cite{Lutz2004dpre} \\  %  Ch.~\ref{ch:pseudo}\\
 \addlinespace
 $(S^2\!\times\!S^2)\# (S^2\!\times\!S^2)$ & $({\mathbb Z},0,{\mathbb Z}^4,0,{\mathbb Z})$ &    $\underline{12}$ & $\underline{12}$ & \cite{KoehlerLutz2004bpre}\\%Ch.~\ref{ch:manifold}\\
 \addlinespace
 ${\mathbb C}{\bf P}^{\,2}\# ({\mathbb C}{\bf P}^{\,2}\# -{\mathbb C}{\bf P}^{\,2})$
                     & $({\mathbb Z},0,{\mathbb Z}^3,0,{\mathbb Z})$ &    13 & ? & \\
 $={\mathbb C}{\bf P}^{\,2}\# (S^2\!\times\!S^2)$ &&&& \\
 \addlinespace
 ${\mathbb C}{\bf P}^{\,2}\#\,(S^2\!\times\!S^2)^{\# 2}$
                     & $({\mathbb Z},0,{\mathbb Z}^5,0,{\mathbb Z})$ &    13 & ? & \\
 \addlinespace
 $(S^3\!\times\!S^1)\#\,{\mathbb C}{\bf P}^{\,2}$
                     & $({\mathbb Z},{\mathbb Z},{\mathbb Z},{\mathbb Z},{\mathbb Z})$ &    14 & ? & \\
 \addlinespace
 $(S^3\hbox{$\times\hspace{-1.62ex}\_\hspace{-.4ex}\_\hspace{.7ex}$}S^1)\#\,{\mathbb C}{\bf P}^{\,2}$
                     & $({\mathbb Z},{\mathbb Z},{\mathbb Z},{\mathbb Z}_2,0)$ &    14 & ? & \\
 \addlinespace
 $({\mathbb C}{\bf P}^{\,2}\#\,{\mathbb C}{\bf P}^{\,2})^{\# 2}$
                     & $({\mathbb Z},0,{\mathbb Z}^4,0,{\mathbb Z})$ &    14 & ? & \\
 \addlinespace
 $(S^2\!\times\!S^2)^{\# 3}$ & $({\mathbb Z},0,{\mathbb Z}^6,0,{\mathbb Z})$ &    14 & ? & \\
 \addlinespace
 $(S^3\!\times\!S^1)\#\,(S^3\!\times\!S^1)$ & $({\mathbb Z},{\mathbb Z}^2,0,{\mathbb Z}^2,{\mathbb Z})$ &    15 & ?  & \\
 \addlinespace
 $(S^3\hbox{$\times\hspace{-1.62ex}\_\hspace{-.4ex}\_\hspace{.7ex}$}S^1)\#\,(S^3\hbox{$\times\hspace{-1.62ex}\_\hspace{-.4ex}\_\hspace{.7ex}$}S^1)$
                     & $({\mathbb Z},{\mathbb Z}^2,0,{\mathbb Z}\oplus{\mathbb Z}_2,0)$ &    15 & ? & \\
 \addlinespace
 $(S^3\hbox{$\times\hspace{-1.62ex}\_\hspace{-.4ex}\_\hspace{.7ex}$}S^1)\#\,({\mathbb C}{\bf P}^{\,2})^{\# 5}$ & $({\mathbb Z},{\mathbb Z},{\mathbb Z}^5,{\mathbb Z}_2,0)$ &    15 & 15 & \cite{KoehlerLutz2004bpre}\\%Ch.~\ref{ch:manifold}\\
 \addlinespace
 $(S^3\!\times\!S^1)\#\,(S^2\!\times\!S^2)$
                     & $({\mathbb Z},{\mathbb Z},{\mathbb Z}^2,{\mathbb Z},{\mathbb Z})$ &    16 & ? & \\
 \addlinespace
 $(S^3\hbox{$\times\hspace{-1.62ex}\_\hspace{-.4ex}\_\hspace{.7ex}$}S^1)\#\,(S^2\!\times\!S^2)$
                     & $({\mathbb Z},{\mathbb Z},{\mathbb Z}^2,{\mathbb Z}_2,0)$ &    16 & ? & \\
 \addlinespace
 ${\mathbb R}{\bf P}^{\,4}$
                     & $({\mathbb Z},{\mathbb Z}_2,0,{\mathbb Z}_2,0)$ &    $\underline{16}$ & ? &  \cite{Lutz2004dpre} \\ %Ch.~\ref{ch:pseudo} \\
 \addlinespace
 K3 surface          & $({\mathbb Z},0,{\mathbb Z}^{22},0,{\mathbb Z})$ &    $\underline{16}$ & $\underline{16}$ & \cite{CasellaKuehnel2001} \\
 \addlinespace
 \addlinespace
\bottomrule
\end{tabular*}
\end{table}

K\"uhnel's bound of Theorem~\ref{thm:kuhnel-4dim} states that $\binom{n-4}{3}\geq 10\,(\chi (M)-2)$
for every combinatorial $4$-manifold $M$ with $n$ vertices (with equality if and only if $M$ is $3$-neighborly).

As a consequence, PL $4$-manifolds $M$ of Euler characteristic $\chi(M)=4$, which include the manifolds
$S^2\!\times\!S^2$, ${\mathbb C}{\bf P}^{\,2}\#\,{\mathbb C}{\bf P}^{\,2}$,
and ${\mathbb C}{\bf P}^{\,2}\# -{\mathbb C}{\bf P}^{\,2}$ with
homo\-logy $H_*=({\mathbb Z},0,{\mathbb Z}^2,0,{\mathbb Z})$,
need at least $10$ vertices for a triangulation.
Since $\binom{10-4}{3}= 10\,(4-2)$ for $n=10$ and $\chi(M)=4$,
every $10$-vertex triangulation of such a manifold would
necessarily be $3$-neighborly according to Theorem~\ref{thm:kuhnel-4dim}.
However, K\"uhnel and Lassmann \cite{KuehnelLassmann1983-unique} proved 
that the boundary of the $5$-simplex and ${\mathbb C}{\bf P}^{\,2}_9$
are the only $3$-neighborly combinatorial $4$-manifolds with $n\leq 13$ vertices.
Thus, at least $11$ vertices are needed for a triangulation of 
a PL $4$-manifold of Euler characteristic $\chi(M)=4$.

%\begin{thm}
%There is a vertex-minimal triangulation of\, $S^2\!\times\!S^2$ with $11$ vertices.\index{manifold!aaS2xS2@$S^2\times S^2$}
%\end{thm}

%\pagebreak

\begin{thm}
The product $S^2\!\times\!S^2$ can be triangulated 
vertex-minimally with $f=(11,\underline{55},150,170,68)$.
\end{thm}

\noindent
\textbf{Proof.} 
%\begin{proof}
We applied bistellar flips to the product triangulation of $S^2\!\times\!S^2$
on $4\cdot 4=16$ vertices and obtained an $11$-vertex triangulation $(S^2\!\times\!S^2)_{11}$
of $S^2\!\times\!S^2$ with $f$-vector $(11,\underline{55},150,170,68)$ and facets:
\begin{center}\small
\begin{tabular}{l@{\hspace{4mm}}l@{\hspace{4mm}}l@{\hspace{4mm}}l@{\hspace{4mm}}l@{\hspace{4mm}}l@{\hspace{4mm}}l@{\hspace{4mm}}l}
%  $12346$       & $12347$       & $12369$       & $12379$       & $12458$       & $12459$       & $12468$ \\
%  $12479$       & $12568$       & $12569$       & $13467$       & $13567$       & $13569$       & $1357\,10$ \\
%  $1359\,11$    & $135\,10\,11$ & $1379\,10$    & $139\,10\,11$ & $1458\,10$    & $1459\,11$    & $145\,10\,11$ \\
%  $1467\,11$    & $1468\,10$    & $146\,10\,11$ & $1479\,11$    & $15678$       & $1578\,10$    &  $1678\,11$ \\
%  $168\,10\,11$ & $178\,10\,11$ & $179\,10\,11$ & $23468$       & $23478$       & $2357\,10$    & $2357\,11$ \\
%  $235\,10\,11$ & $2368\,10$    & $2369\,10$    & $2378\,11$    & $2379\,10$    & $238\,10\,11$ & $24589$ \\
%  $24789$       & $2568\,11$    & $2569\,10$    & $256\,10\,11$ & $25789$       & $2578\,11$    & $2579\,10$ \\
%  $268\,10\,11$ & $3467\,11$    & $3468\,10$    & $3469\,10$    & $3469\,11$    & $3478\,11$    & $3489\,10$ \\
%  $3489\,11$    & $3567\,11$    & $3569\,11$    & $389\,10\,11$ & $4569\,10$    & $4569\,11$    & $456\,10\,11$ \\
%  $4589\,10$    & $4789\,11$    & $5678\,11$    & $5789\,10$    & $789\,10\,11$.
  $12346$ & $12347$ & $12369$ & $12379$ & $12458$ & $12459$ & $12468$ & $12479$ \\
  $12568$ & $12569$ & $13467$ & $13567$ & $13569$ & $1357a$ & $1359b$ & $135ab$ \\
  $1379a$ & $139ab$ & $1458a$ & $1459b$ & $145ab$ & $1467b$ & $1468a$ & $146ab$ \\
  $1479b$ & $15678$ & $1578a$ & $1678b$ & $168ab$ & $178ab$ & $179ab$ & $23468$ \\
  $23478$ & $2357a$ & $2357b$ & $235ab$ & $2368a$ & $2369a$ & $2378b$ & $2379a$ \\
  $238ab$ & $24589$ & $24789$ & $2568b$ & $2569a$ & $256ab$ & $25789$ & $2578b$ \\
  $2579a$ & $268ab$ & $3467b$ & $3468a$ & $3469a$ & $3469b$ & $3478b$ & $3489a$ \\
  $3489b$ & $3567b$ & $3569b$ & $389ab$ & $4569a$ & $4569b$ & $456ab$ & $4589a$ \\
  $4789b$ & $5678b$ & $5789a$ & $789ab$.
\end{tabular}
\end{center}
This triangulation is vertex-minimal by Theorem~\ref{thm:kuhnel-4dim}
and by the result of K\"uhnel and Lassmann \cite{KuehnelLassmann1983-unique}
that there is no $3$-neighborly triangulation of \linebreak
$S^2\!\times\!S^2$ with $10$ vertices.
%\end{proof}
\hfill $\Box$

\begin{conj}
The $f$-vector $(11,\underline{55},150,170,68)$
is com\-po\-nent-wise minimal for
combinatorial triangulations of $S^2\!\times\!S^2$.
\end{conj}

Vertex-transitive triangulations of $S^2\!\times\!S^2$
with $12$ vertices were first found by Sparla \cite{Sparla1997}, \cite{Sparla1998} 
and Lassmann and Sparla~\cite{LassmannSparla2000}.
Altogether, there are three such triangulations; see 
%Chapter~\ref{ch:manifold}. 
\cite{KoehlerLutz2004bpre}.
We applied bistellar flips to these as well and obtained
further minimal triangulations of $S^2\!\times\!S^2$ with $11$ vertices,
which are combinatorially distinct from the above example.
All the examples that we found with $11$ vertices are not symmetric.

\medskip

Combinatorial triangulations of ${\mathbb C}{\bf P}^{\,2}\#\,{\mathbb C}{\bf P}^{\,2}$
and\, ${\mathbb C}{\bf P}^{\,2}\# -{\mathbb C}{\bf P}^{\,2}$
with $9+9-(4+1)=13$ vertices can easily be obtained from K\"uhnel's 
$9$-vertex triangulation ${\mathbb C}{\bf P}^{\,2}_9$ of ${\mathbb C}{\bf P}^{\,2}$
by taking two (disjoint) copies of ${\mathbb C}{\bf P}^{\,2}_9$, removing a
$4$-simplex (with $4+1$ vertices) each, and then gluing both parts
together. Two combinatorially distinct triangulations of ${\mathbb C}{\bf P}^{\,2}\# -{\mathbb C}{\bf P}^{\,2}$
with only $12$ vertices occur as vertex-links of two simply connected $5$-dimensional combinatorial 
pseudomanifolds with homology $H_*=({\mathbb Z},0,0,{\mathbb Z}^{13},0,{\mathbb Z})$; see 
%Chapter~\ref{ch:pseudo}.
\cite{Lutz2004dpre}.
With bistellar flips we obtained triangulations of\, ${\mathbb C}{\bf P}^{\,2}\#\,{\mathbb C}{\bf P}^{\,2}$
and\, ${\mathbb C}{\bf P}^{\,2}\# -{\mathbb C}{\bf P}^{\,2}$ with smaller $f$-vector. 
The resulting lists of facets can be found online at \cite{Lutz_PAGE}.

\begin{prop}
The $4$-manifolds\, ${\mathbb C}{\bf P}^{\,2}\#\,{\mathbb C}{\bf P}^{\,2}$
and\, ${\mathbb C}{\bf P}^{\,2}\# -{\mathbb C}{\bf P}^{\,2}$ can be triangulated 
with $f=(12,57,148,165,66)$.
\end{prop}

\begin{conj}
The $f$-vector $(12,57,148,165,66)$ is com\-po\-nent-wise minimal for
combinatorial triangulations of\, ${\mathbb C}{\bf P}^{\,2}\#\,{\mathbb C}{\bf P}^{\,2}$
and\, ${\mathbb C}{\bf P}^{\,2}\# -{\mathbb C}{\bf P}^{\,2}$.
\end{conj}

We formed further connected sums of the $4$-manifolds ${\mathbb C}{\bf P}^{\,2}$, $S^3\!\times\!S^1$,
$S^3\hbox{$\times\hspace{-1.62ex}\_\hspace{-.4ex}\_\hspace{.7ex}$}S^1$, 
and $S^2\!\times\!S^2$, applied bistellar flips to these, and obtained
small triangulations as listed in Table~\ref{tbl:4-list}.

By enumeration, we also found two vertex-transitive triangulations
$\mbox{}^4\hspace{.3pt}12^{\,2}_{\,1}$ and $\mbox{}^4\hspace{.3pt}12^{\,2}_{\,2}$
of $(S^2\!\times\!S^2)\# (S^2\!\times\!S^2)$ with $12$ vertices; see 
%Chapter~\ref{ch:manifold}. 
\cite{KoehlerLutz2004bpre}.
These triangulations are vertex-minimal according to Theorem~\ref{thm:kuhnel-4dim}.
Moreover, a \emph{tight} (in the sense of \cite{Kuehnel1995-book}) vertex-transitive $15$-vertex 
triangulation $\mbox{}^4\hspace{.3pt}15^{\,4}_{\,1}$ of
$(S^3\hbox{$\times\hspace{-1.62ex}\_\hspace{-.4ex}\_\hspace{.7ex}$}S^1)\#\,({\mathbb C}{\bf P}^{\,2})^{\# 5}$
was obtained; see 
%Chapter~\ref{ch:manifold}. 
\cite{KoehlerLutz2004bpre}.
We believe that this triangulation is vertex-minimal (cf.\ Conjecture~\ref{conj:mintri}).

A tight, vertex-minimal, vertex-transitive triangulation of the K3 surface\index{manifold!aaK3 surface@K3 surface}
was found by Casella and K\"uhnel \cite{CasellaKuehnel2001}. For a survey on the known examples of
tight triangulations see \cite{KuehnelLutz1999}.

%\pagebreak

\section{5-Manifolds}
\label{sec:5-manifolds}

Small triangulations of higher-dimensional manifolds are still rare.
Apart from~$S^d$, triangulated as the boundary of the $(d+1)$-simplex,
the K\"uhnel series \cite{Kuehnel1986a-series},
which contributes a vertex-minimal, vertex-transitive triangulation of the (twisted) 
$S^{d-1}$-bundle over $S^1$ in every dimension $d$, and
the more general series $M^d_k(n)$ of K\"uhnel and Lassmann \cite{KuehnelLassmann1996-bundle}
for $1\leq k\leq d-1$ and $n\geq 2^{d-k}(k+3)-1$,
we know of additional small triangulations only in dimensions\, $d\leq 8$.
The $5$-dimensional examples are listed in Table~\ref{tbl:5-list}.

\begin{table}
\small\centering
\defaultaddspace=0.2em
\caption{\protect\parbox[t]{9.5cm}{Combinatorial $5$-manifolds with
    $n\leq 16$ vertices and smallest known transitive triangulations with $n_{vt}$
    vertices (minimal if underlined).}}\label{tbl:5-list}
\begin{tabular*}{\linewidth}{@{\extracolsep{\fill}}llrrl@{}}
 \addlinespace
 \addlinespace
 \addlinespace
\toprule
 \addlinespace
 \addlinespace
   Manifold          & Homology         &        $n$       &      $n_{vt}$    &  Reference \\ \midrule
 \addlinespace
 \addlinespace
 \addlinespace
     $S^5$           & $({\mathbb Z},0,0,0,0,{\mathbb Z})$ &  $\underline{7}$ &  $\underline{7}$ & \\
 \addlinespace
 $S^3\!\times\!S^2$& $({\mathbb Z},0,{\mathbb Z},{\mathbb Z},0,{\mathbb Z})$ &$\underline{12}$& $\underline{14}$ &  \cite{KoehlerLutz2004bpre}, \cite{Lutz2004apre}\\%Ch.~\ref{ch:manifold}, Ch.~\ref{ch:centrally}\\
 \addlinespace
 $SU(3)/SO(3)$            & $({\mathbb Z},0,{\mathbb Z}_2,0,0,{\mathbb Z})$ &               13 & $\underline{13}$ & \cite{KoehlerLutz2004bpre} \\%Ch.~\ref{ch:manifold}\\
 \addlinespace
 $S^4\hbox{$\times\hspace{-1.62ex}\_\hspace{-.4ex}\_\hspace{.7ex}$}S^1$
                     & $({\mathbb Z},{\mathbb Z},0,0,{\mathbb Z}_2,0)$ & $\underline{13}$ & $\underline{13}$ & \cite{Kuehnel1986a-series} \\
 \addlinespace
 $S^4\!\times\!S^1$  & $({\mathbb Z},{\mathbb Z},0,0,{\mathbb Z},{\mathbb Z})$ &               14 & $\underline{14}$ & \cite{KuehnelLassmann1996-bundle} \\
 \addlinespace
 \addlinespace
\bottomrule
\end{tabular*}
\end{table}

By Theorem~\ref{thm:K-series} we have minimality
for $S^4\hbox{$\times\hspace{-1.62ex}\_\hspace{-.4ex}\_\hspace{.7ex}$}S^1$
with $13$ vertices. For the sphere product $S^4\!\times\!S^1$ it is conjectured 
(cf.\ Conjecture~\ref{conj:products_even_dim}) that $14$ vertices is
best possible.

Vertex-transitive triangulations of the simply connected sphere product $S^3\!\times\!S^2$ with $14$ vertices
and of the simply connected homogeneous $5$-manifold $SU(3)/SO(3)$ with $13$ vertices
were obtained by enumeration; see 
%Chapter~\ref{ch:manifold}. 
\cite{KoehlerLutz2004bpre}.
The unique vertex-transitive triangulation $\mbox{}^5\hspace{.3pt}13^{\,3}_{\,2}$
of $SU(3)/SO(3)$ with $13$ vertices is $3$-neighborly and tight (cf.\ \cite{KuehnelLutz1999}).
Its vertex-minimality is conjectured in 
%Chapter~\ref{ch:manifold}.
\cite{KoehlerLutz2004bpre}.

There are at least four combinatorially distinct vertex-transitive triangulations
of $S^3\!\times\!S^2$ with $14$ vertices. By running the program BISTELLAR
on the example $\mbox{}^5\hspace{.3pt}14^{\,3}_{\,9}$ of 
%Chapter~\ref{ch:manifold}
\cite{KoehlerLutz2004bpre}
we obtained a $12$-vertex triangulation $(S^3\!\times\!S^2)^a_{12}$ of $S^3\!\times\!S^2$ 
with $f$-vector $(12,\underline{66},\underline{220},390,336,112)$ and facets: 
\begin{center}
\small
\begin{longtable}{l@{\hspace{4mm}}l@{\hspace{4mm}}l@{\hspace{4mm}}l@{\hspace{4mm}}l@{\hspace{4mm}}l@{\hspace{4mm}}l@{\hspace{4mm}}l}
  $12346a$ & $12346b$ & $123478$ & $12347b$ & $12348a$ & $12357b$ & $12357c$ & $12359b$ \\
  $12359c$ & $1236ab$ & $12378c$ & $1238ac$ & $1239ab$ & $1239ac$ & $124678$ & $12467b$ \\
  $124689$ & $12469a$ & $12489a$ & $1257bc$ & $1259bc$ & $12678c$ & $1267bc$ & $12689c$ \\
  $1269ab$ & $1269bc$ & $1289ac$ & $134678$ & $13467b$ & $13468a$ & $13579b$ & $13579c$ \\
  $13678c$ & $1367bc$ & $1368ac$ & $136abc$ & $1379ab$ & $1379ac$ & $137abc$ & $145689$ \\
  $14568a$ & $14569a$ & $14589c$ & $1458ac$ & $1459ac$ & $1489ac$ & $15689b$ & $1568ab$ \\
  $1569ab$ & $1579ab$ & $1579ac$ & $157abc$ & $1589bc$ & $158abc$ & $1689bc$ & $168abc$ \\
  $23456a$ & $23456c$ & $23458a$ & $23458b$ & $2345bc$ & $2346bc$ & $23478b$ & $235678$ \\ 
  $23567c$ & $2356a0$ & $23578b$ & $2359bc$ & $23678c$ & $2368ac$ & $236abc$ & $239abc$ \\
  $24567a$ & $24567c$ & $24578a$ & $24578b$ & $2457bc$ & $246789$ & $24679a$ & $2467bc$ \\
  $24789a$ & $25678a$ & $26789a$ & $2689ac$ & $269abc$ & $345679$ & $34567c$ & $345689$ \\
  $34568a$ & $34579c$ & $34589b$ & $3459bc$ & $346789$ & $3467bc$ & $34789b$ & $3479bc$ \\
  $356789$ & $35789b$ & $379abc$ & $45679a$ & $4578ab$ & $4579ac$ & $457abc$ & $4589bc$ \\
  $458abc$ & $4789ab$ & $479abc$ & $489abc$ & $56789a$ & $5689ab$ & $5789ab$ & $689abc$.
\end{longtable}
\end{center}
\vspace{-5mm}
\begin{thm}
The $12$-vertex triangulation $(S^3\!\times\!S^2)^{a}_{12}$ of\, $S^3\times S^2$\, \index{manifold!aaS3xS2@$S^3\times S^2$}
has the minimal number of vertices that a combinatorial
$5$-manifold, different from~$S^5$, 
can have by the Brehm-K\"uhnel bound of Theorem~\ref{thm:bk-minimal}.
In particular, the Brehm-K\"uhnel lower bound is sharp in dimension $5$.
\end{thm}

%\pagebreak

In addition to the example $(S^3\!\times\!S^2)^a_{12}$, we found 
a second triangulation $(S^3\!\times\!S^2)^{b}_{12}$ 
of $S^3\times S^2$ with the same $f$-vec\-tor $(12,\underline{66},\underline{220},390,336,112)$
by applying bistellar flips to the product triangulation of $S^3\!\times\!S^2$
with $20$ vertices. The two triangulations are combinatorially 
distinct: $(S^3\!\times\!S^2)^a_{12}$ has Altshuler-Steinberg determinant  $4471184572226676864$,
whereas the second example $(S^3\!\times\!S^2)^{b}_{12}$ has determinant $4508595451809050112$.
For a list of facets of the second example see 
%Chapter~\ref{ch:tight}.
\cite{KuehnelLutz1999}.
Both triangulations have no non-trivial symmetries: the Altshuler-Steinberg determinants 
of their $12$ vertex links are pairwise distinct, respectively.

Three further minimal triangulations of $S^3\times S^2$
with $12$ vertices were found by starting with 
the vertex-transitive $14$-vertex triangulations $\mbox{}^5\hspace{.3pt}14^{\,3}_{\,13}$,
$\mbox{}^5\hspace{.3pt}14^{\,3}_{\,14}$, and
$\mbox{}^5\hspace{.3pt}14^{\,3}_{\,15}$ from 
%Chapter~\ref{ch:manifold}.
\cite{KoehlerLutz2004bpre}.

With a new and much faster implementation by Nikolaus Witte of the bistellar flip program 
(accessible via the \texttt{TOPAZ} module of the \texttt{polymake} system \cite{polymake})
another twenty examples were obtained by starting from the product triangulation
of $S^3\times S^2$. In fact, we started twenty-six times and each
time achieved a minimal triangulation, but six of these examples 
appeared twice (up to relabeling the vertices). 

\begin{prop}
There are at least\, $25$ combinatorially distinct minimal 
triangulations of\, $S^3\times S^2$ with $12$ vertices.
\end{prop}
The $25$ examples have different Altshuler-Steinberg determinants,
their lists of facets can be found online at \cite{Lutz_PAGE}.

\section{6-Manifolds}
\label{sec:6-manifolds}

There are at least nine combinatorially different vertex-transitive 
$15$-vertex triangulations of $S^3\!\times\!S^3$; see 
%Chapter~\ref{ch:manifold}.
\cite{KoehlerLutz2004bpre}.
Via bistellar flips (by starting with the triangulation $\mbox{}^6\hspace{.3pt}15^{\,7}_{\,2}$
 of $S^3\!\times\!S^3$ from 
%Chapter~\ref{ch:manifold}) 
\cite{KoehlerLutz2004bpre})
we obtained a $13$-vertex triangulation $(S^3\!\times\!S^3)^{a}_{13}$
of $S^3\!\times\!S^3$ with $f$-vector $(13,\underline{78},\underline{286},\underline{715},1014,728,208)$
and facets:
%\vspace{-1mm}
\begin{center}\small
\begin{longtable}{@{}l@{\hspace{3mm}}l@{\hspace{3mm}}l@{\hspace{3mm}}l@{\hspace{3mm}}l@{\hspace{3mm}}l@{\hspace{3mm}}l@{\hspace{3mm}}l@{}}
$123456c$\! & $123456d$\! & $12345ab$\! & $12345ac$\! & $12345bd$\! & $12346cd$\! & $123479a$\! & $123479d$ \\
$12347ad$\! & $12348ab$\! & $12348ad$\! & $12348bd$\! & $12349ac$\! & $12349cd$\! & $12356bc$\! & $12356bd$ \\
$1235abc$\! & $1236bcd$\! & $12379ad$\! & $12389ab$\! & $12389ad$\! & $12389bd$\! & $1239abc$\! & $1239bcd$ \\
$1245678$\! & $124567b$\! & $1245689$\! & $124569a$\! & $12456ac$\! & $12456bd$\! & $1245789$\! & $124579a$ \\
$12457ab$\! & $124678b$\! & $124689b$\! & $12469ac$\! & $12469bd$\! & $12469cd$\! & $124789d$\! & $12478ab$ \\
$12478ad$\! & $12489bd$\! & $125678c$\! & $12567bc$\! & $125689a$\! & $12568ac$\! & $125789a$\! & $12578ac$ \\
$1257abc$\! & $12678ab$\! & $12678ac$\! & $1267abc$\! & $12689ab$\! & $1269abc$\! & $1269bcd$\! & $12789ad$ \\
$1345679$\! & $134567d$\! & $134569a$\! & $13456ac$\! & $134579a$\! & $13457ad$\! & $1345abd$\! & $134679c$ \\
$13467cd$\! & $13469ac$\! & $13479cd$\! & $1348abd$\! & $1356789$\! & $135678c$\! & $13567bc$\! & $13567bd$ \\
$135689a$\! & $13568ac$\! & $135789c$\! & $13579ad$\! & $13579bc$\! & $13579bd$\! & $13589ad$\! & $13589bc$ \\
$13589bd$\! & $1358abc$\! & $1358abd$\! & $136789c$\! & $1367bcd$\! & $13689ac$\! & $1379bcd$\! & $1389abc$ \\
$1456789$\! & $14567bd$\! & $1457abd$\! & $146789c$\! & $14678bd$\! & $14678cd$\! & $14689bc$\! & $1468bcd$ \\
$1469bcd$\! & $14789cd$\! & $1478abd$\! & $1489bcd$\! & $15789ad$\! & $15789cd$\! & $1578acd$\! & $1579bcd$ \\
$157abcd$\! & $1589bcd$\! & $158abcd$\! & $1678abd$\! & $1678acd$\! & $167abcd$\! & $1689abc$\! & $168abcd$ \\
$23456cd$\! & $23458bc$\! & $23458bd$\! & $23458cd$\! & $2345abc$\! & $23478ab$\! & $23478ad$\! & $23478bc$ \\
$23478cd$\! & $23479ac$\! & $23479cd$\! & $2347abc$\! & $235678b$\! & $235678c$\! & $23567bc$\! & $235689b$ \\
$235689d$\! & $23568cd$\! & $23569bd$\! & $23578bc$\! & $23589bd$\! & $23678ab$\! & $23678ad$\! & $23678cd$ \\
$23679ab$\! & $23679ad$\! & $23679bd$\! & $2367bcd$\! & $23689ab$\! & $23689ad$\! & $2379abc$\! & $2379bcd$ \\
$245678b$\! & $245689b$\! & $24569ac$\! & $24569bd$\! & $24569cd$\! & $245789c$\! & $24578bc$\! & $24579ac$ \\
$2457abc$\! & $24589bd$\! & $24589cd$\! & $24789cd$\! & $25689ad$\! & $2568acd$\! & $2569acd$\! & $25789ad$ \\
$25789cd$\! & $2578acd$\! & $2579acd$\! & $2678acd$\! & $2679abd$\! & $267abcd$\! & $269abcd$\! & $279abcd$ \\
$345679a$\! & $34567ad$\! & $3456acd$\! & $3458abc$\! & $3458abd$\! & $3458acd$\! & $346789b$\! & $346789c$ \\
$34678ab$\! & $34678ad$\! & $34678cd$\! & $34679ab$\! & $34689ab$\! & $34689ac$\! & $3468acd$\! & $34789bc$ \\
$3479abc$\! & $3489abc$\! & $356789b$\! & $35679ad$\! & $35679bd$\! & $35689ad$\! & $3568acd$\! & $35789bc$ \\
$456789b$\! & $45679ab$\! & $4567abd$\! & $4569abd$\! & $4569acd$\! & $45789bc$\! & $4579abc$\! & $4589bcd$ \\
$458abcd$\! & $459abcd$\! & $4678abd$\! & $4689abc$\! & $468abcd$\! & $469abcd$\! & $5679abd$\! & $579abcd$.
\end{longtable}
\end{center}
\vspace{-5mm}
\begin{thm}
The $13$-vertex triangulation $(S^3\!\times\!S^3)^{a}_{13}$ of\, $S^3\times S^3$\, \index{manifold!aaS3xS3@$S^3\times S^3$}
has the minimal number of vertices that a combinatorial
$6$-manifold, different from~$S^6$, can have by the Brehm-K\"uhnel bound of Theorem~\ref{thm:bk-minimal}.
In particular, the Brehm-K\"uhnel lower bound is sharp in dimension $6$.
\end{thm}

Another $13$-vertex triangulation $(S^3\!\times\!S^3)^{b}_{13}$ of $S^3\!\times\!S^3$ was obtained
by starting the bistellar flip program on the vertex-transitive triangulation 
$\mbox{}^6\hspace{.3pt}15^{\,7}_{\,1}$ of 
%Chapter~\ref{ch:manifold}.
\cite{KoehlerLutz2004bpre}.
Both $13$-vertex triangulations of\, $S^3\!\times\!S^3$
are $4$-neighborly and thus are tight, since equality holds in Proposition~4.6 of \cite{Kuehnel1995-book};
see also \cite{KuehnelLutz1999}. The facets of $(S^3\!\times\!S^3)^{b}_{13}$ 
can be found in 
%Chapter~\ref{ch:tight}.
\cite{KuehnelLutz1999}.
The two examples $(S^3\!\times\!S^3)^{a}_{13}$ and
$(S^3\!\times\!S^3)^{b}_{13}$ have the same Altshuler-Steinberg determinant $745714154823444619853824$.
However, the Altshuler-Steinberg determinants of their vertex links
are pairwise distinct. It follows that the two examples are asymmetric.

Two further tight $13$-vertex triangulations of $S^3\!\times\!S^3$
were obtained with bistellar flips by starting from the
centrally-symmetric $16$-vertex triangulations\, $\mbox{}^6_{\times}\hspace{.3pt}16^{\,cy}_{\,1}$
and\, $\mbox{}^6_{\times}\hspace{.3pt}16^{\,di}_{\,1}$ of 
%Chapter~\ref{ch:centrally}.
\cite{Lutz2004apre}.
The resulting minimal triangulations have the same Altshuler-Steinberg determinant as the two
examples before, but again, the Altshuler-Steinberg determinants of their vertex links
differ.

\begin{prop}
There are at least\, $4$ combinatorially distinct minimal 
triangulations of\, $S^3\times S^3$ with $13$ vertices.
\end{prop}

The short list of $6$-manifolds for which we know small triangulations
is given in Table~\ref{tbl:6-list}.

\begin{table}
\small\centering
\defaultaddspace=0.2em
\caption{\protect\parbox[t]{8.1cm}{Combinatorial $6$-manifolds with
    $n\leq 16$ vertices and smallest known transitive triangulations with $n_{vt}$
    vertices (minimal if underlined).}}\label{tbl:6-list}
\begin{tabular}{@{}l@{\hspace{10mm}}r@{\hspace{10mm}}r@{\hspace{10mm}}l@{}}
 \addlinespace
 \addlinespace
 \addlinespace
\toprule
 \addlinespace
 \addlinespace
   Manifold          &        $n$       &      $n_{vt}$    &  Reference \\ \midrule
 \addlinespace
 \addlinespace
 \addlinespace
     $S^6$           &  $\underline{8}$ &  $\underline{8}$ & \\
 \addlinespace
 $S^3\!\times\!S^3$  & $\underline{13}$ &               15 & \cite{KoehlerLutz2004bpre} \\%Ch.~\ref{ch:manifold}\\
 \addlinespace
 $S^5\!\times\!S^1$  & $\underline{15}$ & $\underline{15}$ & \cite{Kuehnel1986a-series} \\
 \addlinespace
 \addlinespace
\bottomrule
\end{tabular}
\end{table}

\section{7-Manifolds}
\label{sec:7-manifolds}

Vertex-transitive triangulations of $7$-manifolds 
with up to $15$ vertices are enumerated in 
%Chapter~\ref{ch:manifold}
\cite{KoehlerLutz2004bpre} 
(with the
exception of possible examples corresponding to the transitive actions
of the groups ${\mathbb Z}_{14}$, $D_7$, and ${\mathbb Z}_{15}$ on $14$ and $15$ vertices).
All resulting examples are triangulations of $S^7$. 
(It is open whether there are vertex-transitive
triangulations of $7$-manifolds, different from $S^7$, with a 
cyclic automorphism group on $15$ vertices.)

Triangulations of centrally symmetric $7$-manifolds
with a vertex-transitive cyclic group action on $18$ vertices
are enumerated in 
%Chapter~\ref{ch:centrally}: 
\cite{Lutz2004apre}:
There is one vertex-transitive
centrally symmetric $18$-vertex triangulation with cyclic symmetry of $S^4\!\times\!S^3$
and one of $S^5\!\times\!S^2$ each.
These are the smallest triangulations for these manifolds
that have been achieved so far; see Table~\ref{tbl:7-list}.

\begin{table}
\small\centering
\defaultaddspace=0.2em
\caption{\protect\parbox[t]{8.1cm}{Combinatorial $7$-manifolds with
    $n\leq 20$ vertices and smallest known vertex-transitive triangulations with $n_{vt}$
    vertices (minimal if underlined).}}\label{tbl:7-list}
\begin{tabular}{@{}l@{\hspace{10mm}}r@{\hspace{10mm}}r@{\hspace{10mm}}l@{}}
 \addlinespace
 \addlinespace
 \addlinespace
\toprule
 \addlinespace
 \addlinespace
   Manifold          &        $n$       &      $n_{vt}$    &  Reference \\ \midrule
 \addlinespace
 \addlinespace
 \addlinespace
      $S^7$           &  $\underline{9}$ &  $\underline{9}$ & \\
 \addlinespace
$S^6\hbox{$\times\hspace{-1.62ex}\_\hspace{-.4ex}\_\hspace{.7ex}$}S^1$
                     & $\underline{17}$ & $\underline{17}$ & \cite{Kuehnel1986a-series} \\
 \addlinespace
 $S^6\!\times\!S^1$  &               18 &               18 & \cite{KuehnelLassmann1996-bundle} \\
 \addlinespace
 $S^5\!\times\!S^2$  &               18 &               18 & \cite{Lutz2004apre}\\%Ch.~\ref{ch:centrally} \\
 \addlinespace
 $S^4\!\times\!S^3$  &               18 &               18 &  \cite{Lutz2004apre}\\%Ch.~\ref{ch:centrally} \\
 \addlinespace
 \addlinespace
\bottomrule
\end{tabular}
\end{table}

\section{8-Manifolds}
\label{sec:8-manifolds}

According to the Brehm-K\"uhnel bound of Theorem~\ref{thm:bk-minimal}(a),
a combinatorial $8$-manifold, different from $S^8$, has at least $15$ vertices.\index{projective space!`quaternionic'}
It is a `manifold like the quaternionic projective plane' if\, it has
$15$ vertices. Such an example  $M^8_{15}$ with a vertex-transitive $A_5$-action and 
two further non-transitive examples, which are PL homeomorphic
to the transitive one, were found by Brehm and K\"uhnel \cite{BrehmKuehnel1992}. 
We denote their $8$-manifold by ${\sim}{\mathbb H}{\bf P}^{\,2}$.
With bistellar flips we found three further triangulations of ${\sim}{\mathbb H}{\bf P}^{\,2}$;
see \cite{Lutz_PAGE} for their lists of facets.

\begin{prop}
There are at least\, $6$ combinatorially distinct vertex-min\-i\-mal 
triangulations of\, the Brehm and K\"uhnel manifold 
${\sim}{\mathbb H}{\bf P}^{\,2}$
with $15$ vertices.
\end{prop}

Centrally symmetric $20$-vertex triangulations of\, $S^4\times S^4$
and of\, $S^5\times S^3$  (one each) 
with dihedral symmetry are given in 
%Chapter~\ref{ch:centrally}.
\cite{Lutz2004apre}.
These are the smallest known triangulations for these two manifolds; see Table~\ref{tbl:7-list}.
\begin{table}
\small\centering
\defaultaddspace=0.1em
\caption{\protect\parbox[t]{8.1cm}{Combinatorial $8$-manifolds with
    $n\leq 20$ vertices and smallest known vertex-transitive triangulations with $n_{vt}$
    vertices (minimal if underlined).}}\label{tbl:8-list}
\begin{tabular}{@{}l@{\hspace{10mm}}r@{\hspace{10mm}}r@{\hspace{10mm}}l@{}}
 \addlinespace
 \addlinespace
 \addlinespace
 \addlinespace
 \addlinespace
 \addlinespace
\toprule
 \addlinespace
 \addlinespace
   Manifold          &        $n$       &      $n_{vt}$    &  Reference \\ \midrule
 \addlinespace
 \addlinespace
 \addlinespace
     $S^8$           & $\underline{10}$ & $\underline{10}$ & \\
 \addlinespace
${\sim}{\mathbb H}{\bf P}^{\,2}$
                     & $\underline{15}$ & $\underline{15}$ & \cite{BrehmKuehnel1992} \\
 \addlinespace
 $S^7\!\times\!S^1$  & $\underline{19}$ & $\underline{19}$ & \cite{Kuehnel1986a-series} \\
 \addlinespace
 $S^5\!\times\!S^3$  &             20 &               20 &  \cite{Lutz2004apre}\\%Ch.~\ref{ch:centrally}\\ 
 \addlinespace
 $S^4\!\times\!S^4$  &             20 &               20 &  \cite{Lutz2004apre}\\%Ch.~\ref{ch:centrally}\\ 
 \addlinespace
 \addlinespace
\bottomrule
\end{tabular}
\end{table}

%%%%%%%%%%%%%%%%%%%%%%%%%%%%%%%%%%%%%%%%%%%%%%%%%%%%%%%%%%%%%%%%%%%

\pagebreak

%\bibliography{../../../references}
\bibliography{}

\bigskip
\medskip

\noindent
Frank H.\ Lutz\\
Technische Universit\"at Berlin\\
Fakult\"at II - Mathematik und Naturwissenschaften\\
Institut f\"ur Mathematik, Sekr. MA 6-2\\
Stra\ss e des 17.\ Juni 136\\
D-10623 Berlin\\
{\tt lutz@math.tu-berlin.de}

\end{document}